\documentclass[11pt,reqno]{amsart}
\usepackage{amsmath,amsfonts,amssymb,amsthm}
\usepackage[mathcal]{euscript}
\usepackage{mathrsfs}
\usepackage[tmargin=1.4in,bmargin=1.4in,lmargin=1.2in,rmargin=1.2in]{geometry}
\usepackage{graphicx}
\usepackage{tikz}
\usepackage{cite}
\graphicspath{{./Images/}}
\usepackage[font=small,format=plain,labelfont=sf,up,textfont=sf,up]{caption}
\usepackage{boxedminipage}

\allowdisplaybreaks

\usepackage{microtype}
\usepackage{bbm}
\usepackage{color}
\usepackage[textsize=small]{todonotes}

\usepackage{thmtools}
\declaretheorem[numberwithin=section, name=Theorem]{thm}
\declaretheorem[sibling=thm, name=Lemma]{lem}

\declaretheorem[sibling=thm, name=Definition]{defi}
\declaretheorem[sibling=thm, name=Proposition]{prop}
\declaretheorem[sibling=thm, name=Corollary]{cor}

\definecolor{halfgray}{gray}{0.55}
\definecolor{webgreen}{rgb}{0,.5,0}
\definecolor{webbrown}{rgb}{.6,0,0}
\definecolor{Maroon}{cmyk}{0, 0.87, 0.68, 0.32}
\definecolor{royalblue}{RGB}{0,0,139}
\definecolor{lightblue}{RGB}{150,230,255}
\definecolor{Black}{cmyk}{0, 0, 0, 0}
\definecolor{pinkish}{RGB}{255, 192, 203}

\usepackage{hyperref}
\hypersetup{%
    colorlinks=true, linktocpage=true, pdfstartpage=1, pdfstartview=FitV,%
    breaklinks=true, pdfpagemode=UseNone, pageanchor=true, pdfpagemode=UseOutlines,%
    plainpages=false, bookmarksnumbered, bookmarksopen=true, bookmarksopenlevel=1,%
    hypertexnames=true, pdfhighlight=/O,
    urlcolor=webbrown, linkcolor=royalblue, citecolor=webgreen, 
}

\usepackage{tikz}
\pgfdeclarelayer{nodelayer}
\pgfdeclarelayer{edgelayer}
\pgfsetlayers{edgelayer,nodelayer,main}
\tikzstyle{n}=[circle,fill=black,draw=black,line width=0.2 pt,minimum size=0.1 cm ]

\numberwithin{equation}{section}

\newcommand{\shift}{\!\!\!\!}
\newcommand{\sss}{\scriptscriptstyle}
\newcommand{\Var}{\mathrm{Var}}
\renewcommand{\P}{\mathbb{P}}

\newcommand{\indi}{\mathbbm{1}}

\newcommand{\cluster}{\mathscr{C}}

\newcommand{\E}{\mathbb{E}}
\newcommand{\e}{\mathrm{e}}

\renewcommand{\P}{\mathbb{P}}

\newcommand{\Z}{\mathbb{Z}}

\newcommand{\bbN}{\mathbb{N}}

\newcommand{\vep}{\varepsilon}

\newcommand{\Att}{\mathtt{A}}
\newcommand{\Dtt}{\mathtt{D}}
\newcommand{\Ptt}{\mathtt{P}}
\newcommand{\Fcal}{\mathcal{F}}
\newcommand{\Ecal}{\mathcal{E}}

\newcommand{\Pcal}{\mathcal{P}}
\newcommand{\Ucal}{\mathcal{U}}
\newcommand{\Vcal}{\mathcal{V}}
\newcommand{\Tcal}{\mathcal{T}}

\newcommand{\Hcal}{\mathcal{H}}
\newcommand{\Lcal}{\mathcal{L}}

\newcommand{\Wcal}{\mathcal{W}}

\newcommand{\Cbf}{\mathbf{C}}
\newcommand{\Dbf}{\mathbf{D}}
\newcommand{\Abf}{\mathbf{A}}
\newcommand{\Asf}{\mathsf{A}}
\newcommand{\Dsf}{\mathsf{D}}
\newcommand{\Psf}{\mathsf{P}}
\newcommand{\Nbf}{\mathbf{N}}
\newcommand{\Xbf}{\mathbf{X}}
\newcommand{\Cctv}{\tilde\cluster(v)}

\newcommand{\Bin}{{\sf{Bin}}}

\newcommand{\conn}{\longleftrightarrow}

\newcommand{\tmix}{{t_{\mathrm{mix}}}}

\newcommand{\eqan}[1]{\begin{align}#1\end{align}}
\newcommand{\nn}{\nonumber}

\newcommand{\Pp}{\mathbb{P}_p}
\newcommand{\Ep}{\mathbb{E}_p}
\newcommand{\Zd}{\mathbb{Z}^d}

\newcommand{\chil}{\chi_{\sss \mathrm{line}}}
\newcommand{\Sp}{\mathbf{Sp}}
\newcommand{\Spr}{\mathrm{Sp}}

\def\arrowfillCS#1#2#3#4{
   \thickmuskip0mu\medmuskip\thickmuskip\thinmuskip\thickmuskip
   \relax#4#1\mkern-7mu
   \cleaders\hbox{$#4\mkern-2mu#2\mkern-2mu$}\hfill
   \mkern-7mu#3
}
\def\lrfill{\arrowfillCS\leftarrow\relbar\rightarrow\relax}

\newcommand{\arf}[1]{\stackrel{#1}{\longleftrightarrow}}
\newcommand{\arfl}[1]{\stackrel{\,\,\, #1 \,\,\, }{\lrfill}}

\usepackage[normalem]{ulem}

\newcommand{\ch}[1]{#1}
\newcommand{\st}[1]{#1}
\newcommand{\stm}[1]{}

\begin{document}

\title[Expansion of percolation critical points for Hamming graphs]
{Expansion of percolation critical points \\
for Hamming graphs}

\author[L.\ Federico]{Lorenzo Federico}
\email{l.federico@tue.nl}

\author[R.\ van der Hofstad]{Remco van der Hofstad}
\email{r.w.v.d.hofstad@tue.nl}

\author[F.\ den Hollander]{Frank den Hollander}
\email{denholla@math.leidenuniv.nl}

\author[T.\ Hulshof]{Tim Hulshof}
\email{w.j.t.hulshof@tue.nl}

\begin{abstract}
The Hamming graph $H(d,n)$ is the Cartesian product of $d$ complete graphs on $n$ vertices. 
Let $m=d(n-1)$ be the degree and $V = n^d$ be the number of vertices of $H(d,n)$. Let $p_c^{(d)}$ be the critical point for bond percolation on $H(d,n)$. We show that, for $d \in \bbN$ fixed and $n \to \infty$,
	\begin{equation*}
	p_c^{(d)}= \dfrac{1}{m} + \dfrac{2d^2-1}{2(d-1)^2}\dfrac{1}{m^2} 
	+ O(m^{-3}) + O(m^{-1}V^{-1/3}),
	\end{equation*}
which extends the asymptotics found in \cite{BorChaHofSlaSpe05b} by one order. The term 
$O(m^{-1}V^{-1/3})$ is the width of the critical window. For $d=4,5,6$ we have $m^{-3} = 
O(m^{-1}V^{-1/3})$, and so the above formula represents the full asymptotic expansion of 
$p_c^{(d)}$. In \cite{FedHofHolHul16a} \st{we show that} this formula is a crucial ingredient in 
the study of critical bond percolation on $H(d,n)$ for $d=2,3,4$. The proof uses a lace expansion 
for the upper bound and a novel comparison with a branching random walk for the lower bound. 
The proof of the lower bound also yields a refined asymptotics for the susceptibility of a subcritical 
Erd\H{o}s-R\'enyi random graph. 
\end{abstract}

\date{\today}

\maketitle

\noindent
{\small {\it MSC 2010.} 60K35, 60K37, 82B43.

\noindent
{\it Keywords and phrases.}
Hamming graph, percolation, critical point, critical window, lace expansion.}

\medskip


\section{Introduction and main result}


\subsection{Percolation on the Hamming graph}
The Hamming graph $H(d,n)$ is the Cartesian product of $d$ complete graphs on $n$ vertices 
(e.g., $H(3,7) = K_7 \times K_7 \times K_7$). Bernoulli bond percolation is the model where, 
given a graph, each edge is retained independently with the same probability $p$. In this paper 
we study the location of the critical point \ch{of bond percolation on $H(d,n)$} for the phase transition in the size of the largest connected 
component when $d$ is fixed and $n \to \infty$. 

Formally, we define the Hamming graph $H(d,n)$ for $d,n \in \bbN$ as the graph with vertex 
set $\Vcal := \{0,1,\dots,n-1\}^d$ and edge set
\begin{equation}
	\Ecal := \{ (v,w) : v,w \in \Vcal, \ v_j \neq w_j \text{ for exactly one }j \}. 
\end{equation}
Thus, $H(d,n)$ is a transitive graph on $V:= n^d$ vertices with degree $m:=d(n-1)$. Bernoulli 
bond percolation is synonymous with the probability space $(\Omega, \Pp)$, where $\Omega 
:= \{0,1\}^{\Ecal}$ and $\Pp$ is the measure such that
\begin{equation}
	\Pp(\omega) = \prod_{e \in \Ecal} \big((1-p) \delta_{0,\omega(e)} 
	+  p \delta_{1,\omega(e)} \big)
	 \qquad \forall\,\omega \in \Omega,
\end{equation}
where $\delta_{x,y}$ is the Kronecker delta. When $\omega(e)=1$ we say that the edge $e$ is 
\emph{open}, when $\omega(e)=0$ we say that the edge $e$ is \emph{closed}. Given a vertex 
$x \in \Vcal$, we write $\cluster(x)$ for the \ch{graph whose vertex set consists of all vertices that can be reached from $x$ through 
a path of open edges, and whose edge set consists of all open edges between these vertices.} 
We call $\cluster(x)$ the \emph{connected component} of $x$, or \emph{cluster} 
of $x$, \ch{and write $|\cluster(x)|$ for its number of vertices.} We write $\cluster_1$ for the cluster \ch{$\cluster(x)$} with the largest cardinality \ch{$|\cluster(x)|$} (using some tie-breaking 
rule). Two of the main objects of study in percolation are $|\cluster(x)|$ and $|\cluster_1|$, the cardinalities 
of $\cluster(x)$ and $\cluster_1$. For percolation on \emph{infinite graphs} $G$ it is often observed 
that the \emph{critical point of the percolation phase transition on $G$}, defined by
\begin{equation}
	p_c^{G} := \inf\{ p \in [0,1]  \colon  \Pp(|\cluster\ch{(x)}| =\infty) >0\},
\end{equation}
is non-trivial, i.e., $p_c^G \in (0,1)$ (see for example Grimmett \cite{Grim12}) for most infinite graphs (an 
exception being $\mathbb{Z}^1$). Moreover,  Aizenman and Barsky \cite{AizBar87} and independently Menshikov 
\cite{Mens86} proved  that on transitive 
graphs, 
\begin{equation}
	p_c^{G} = \sup\{ p \in [0,1] \colon  \Ep[|\cluster(x)|] < \infty \}.
\end{equation}

Since we consider percolation on $H(d,n)$ with $d,n$ finite and $\Pp$ is a product measure, 
any event that is measurable with respect to $\Pp$ has a probability that is a polynomial in 
$p$, and therefore is continuous in $p$: the finite model cannot undergo a non-trivial phase 
transition in $p$ as described above. Nevertheless, it \emph{does} make sense to study the 
percolation phase transition on finite graphs in the limit as $n\to\infty$. To see why, let us 
give a rough sketch of an important related problem: the emergence of the giant component 
in the \emph{Erd\H{o}s-R\'enyi Random Graph} (ERRG).


\subsection{Giant component}

The Erd\H{o}s-R\'enyi random graph is the common name for percolation on the complete 
graph $K_n$. Erd\H{o}s and R\'enyi \cite{ErdRen59} proved that in the limit as $n \to \infty$, 
if $p = p(n) < n^{-1}$, then $|\cluster_1| = \Theta(\log n)$ w.h.p.,
\footnote{Given a sequence of random variables $(X_n)_{n \in 
\bbN}$, we write $X_n = \Theta(f(n))$ w.h.p.\ (with high probability) if there exist constants 
$C \ge c > 0$ such that $\Pp(c f(n) \le X_n \le C f(n)) \to 1$ as $n \to \infty$.}
while if $p > n^{-1}$, then $|\cluster _1| = \Theta(n)$ w.h.p. Moreover, zooming in on the transition point $n^{-1}$ by choosing $p = (1+ \vep_n)n^{-1}$ for a sequence $(\vep_n)_{n\in\bbN}$ such that $\lim_{n\to\infty} \vep_n = 0$, Bollob\'as \cite{Bol84} showed that 
\footnote{Subsequent results in \cite{LucPitWie94,Ald97,Pit01,NacPer07} are much sharper and comprehensive than what is summarized here, and there is an extensive body of literature on the problem.}
\begin{itemize}
	\item $|\cluster_1| = \Theta(\vep_n^{-2} \log(\vep_n^3 n))$ w.h.p.\ when 
		$\vep_n^3 n \to - \infty$ (subcritical),
	\item $|\cluster_1| = \Theta(n^{2/3})$ w.h.p.\ when $\vep_n n^3 \to a \in \mathbb{R}$
	(critical),
	\item $|\cluster_1| = \Theta(\vep_n n)$ w.h.p.\ when $\vep_n n^3 \to +\infty$ (supercritical).
\end{itemize}
What this shows is that the size of the largest component undergoes a sharp transition around~$n^{-1}$. As mentioned above, there is no critical point for a finite graph, but the transition occurs 
in a slice of the parameter space with a width of order $n^{-4/3}$, which is asymptotically vanishing 
with respect to the center of the window located around $n^{-1}$. This behaviour inspired the 
notion of \emph{critical window}: to indicate that the transition of the ERRG occurs around $n^{-1}$ 
in a range of width $n^{-4/3}$, we use the short-hand notation
\footnote{Given three sequences $(a_n), (b_n), (c_n)$, we write that $a_n = b_n + O(c_n)$ when 
there exists a constant $K < \infty$ such that $|a_n-b_n|\leq K c_n$ for all $n$.}
\begin{equation}
\label{e:critwindowERRG}
	p_c^{K_n} = n^{-1} + O(n^{-4/3}).
\end{equation}

Erd\H{o}s and Spencer \cite{ErdSpe79} conjectured that if we replace $K_n$ by a more ``geometric'' 
graph sequence (their primary candidate was  $H(d,2)$, the $d$-dimensional \emph{hypercube}, 
with $d \to \infty$), then the critical behaviour should remain largely intact. In fact, it turned out that 
to a large extent the picture is the same for a large class of graph sequences with ``sufficiently weak'' 
geometries.

Of particular interest to us here are the papers by Borgs et al.\ \cite{BorChaHofSlaSpe05a , 
BorChaHofSlaSpe05b , BorChaHofSlaSpe06}, demonstrating that 
graph sequences satisfying the so-called \emph{triangle condition} (which serves as an indicator 
of what is meant by sufficiently weak geometry; see e.g.\  \cite{AizNew84 , BorChaHofSlaSpe05a , 
BorChaHofSlaSpe05b , BorChaHofSlaSpe06 , HarSla90}) have a phase transition that strongly 
resembles that of the ERRG, and that both $(H(d,2))_{d \in \bbN}$ and $(H(d,n))_{n \in \bbN}$ 
satisfy the triangle condition. \ch{More precisely, consider a sequence $(G_m)_{m \in\mathbb{N}}
= (\Vcal_m, \Ecal_m)_{m\in\bbN}$ of vertex transitive graphs of degree $m$, and write $V_m 
:= |\Vcal_m|$. Write $x \conn y$ for the event that $y \in \cluster(x)$, and define the \emph{two-point function} $\tau_p(x-y) := \Pp(x \conn y)$ and the \emph{susceptibility} $\chi(p) := \Ep[|\cluster(x)|] =\sum_{y\in\Vcal_m} \tau_p(x-y)$ (note that $\chi(p)$ does not depend on $x$ by transitivity and that $\tau_p(x-y)$ depends on the relative difference of $x$ and $y$ only because the graphs under consideration are tori).}  The triangle condition is satisfied for percolation on $(G_m)$ if for all $p$ such that $\chi(p)^3/V_m \le \beta_0$ for some sufficiently small $\beta_0$, and for all $x,y \in \Vcal_m$, we have
\footnote{Here and below we will frequently suppress sub- and superscripts when their presence is clear from context. Likewise, we do not always stress that we are considering asymptotic results for sequences.}
\begin{equation}
\label{e:triangledef}
	\nabla_p(x,y) := \sum_{u,v} \tau_p(x -u) \tau_p(u-v) \tau_p(v - y) 
	= \delta_{x,y} + 10\,\frac{\chi(p)^3}{V} + O(m^{-1}).
\end{equation}
Borgs et al.\ prove that the triangle condition holds for a class of models that includes 
$(H(d,n))_{n \in \bbN}$ for any fixed $d \ge 2$ (see \cite[Theorem 1.3]{BorChaHofSlaSpe05b}). 
\ch{An alternative proof, applying to e.g.\ Hamming graphs and hypercubes, was given by van der Hofstad and Nachmias \cite{HofNac12,HofNac13}.}


\subsection{Critical window}
 
Fix some $\theta \in (0, \infty)$ and define $ p_c^{G_m}(\theta)$ as the unique solution of
the equation
\begin{equation}
\label{e:pcdef}
	\chi(p_c(\theta)) = \theta V^{1/3}.
\end{equation}
Borgs et al.\ \cite{BorChaHofSlaSpe05a, BorChaHofSlaSpe05b} prove that if we consider 
percolation on a sequence $(G_m)$ that satisfies the triangle condition \eqref{e:triangledef} with 
$p = p_c(\theta) (1 + \vep_m)$ and $\vep_m  \to 0$, then we see subcritical behaviour when 
$\vep^3 V \to - \infty$ and critical behaviour when $\vep^3 V \to a \in \mathbb{R}$, just as in the 
ERRG. Sharper results about mean-field supercritical behaviour of percolation models when 
$\vep^3 V \to \infty$ were derived later by van der Hofstad and Nachmias \cite{HofNac12}, 
\ch{who investigate the supercritical phase and thus establish that \eqref{e:pcdef} really constitutes the critical window
for several high-dimensional tori including the hypercube and Hamming graphs.} 
Moreover, it was shown in \cite[Theorem 1.1]{BorChaHofSlaSpe05a} that the critical window satisfies
\begin{equation}
\label{e:critwindowBCHHS}
	p_c(\theta) = m^{-1} +O(m^{-2}) + O(m^{-1} V^{-1/3}),
\end{equation}
and that $p_c(\theta_1)-p_c(\theta_2)=O(m^{-1}V^{-1/3})$ for any $\theta_1,\theta_2 >0$, i.e.,
any choice of $\theta$ yields \emph{the same} critical window.

Compare \eqref{e:critwindowBCHHS} with the critical window of the ERRG in \eqref{e:critwindowERRG}, 
and note that $K_n$ has $O(m^{-1} V^{-1/3}) = O(n^{-4/3})$ because $m = n-1$ and $V = n$. Thus, 
by that analogy, the second error term above corresponds to the width of the critical window, while the 
first error term can be viewed as a ``correction'' in $m^{-1}$ to $p_c$ itself. In this interpretation, 
\eqref{e:critwindowBCHHS} describes the critical window \emph{asymptotically precisely} for the 
two-dimensional Hamming graph $H(2,n)$, since in this case $m = 2(n-1)$ and $V =n^2$, so that 
the correction term $m^{-2}$ is vanishingly small compared to $m^{-1} V^{-1/3}$. Moreover, 
\eqref{e:critwindowBCHHS} is also asymptotically precise for $H(3,n)$ \ch{because the two $O$-terms coincide.}


\subsection{Expansion of the critical point}

This brings us to the main result of our paper. We write $p_c^{\sss (d)}(\theta)$ for the critical value of percolation on $H(d,n)$ defined in \eqref{e:critwindowBCHHS}, 
and compute the second term of $p_c^{\sss (d)} (\theta)$ for all $d \ge 2$:

\begin{thm}[\ch{Critical window for percolation on $H(d,n)$}]
\label{thm-pc} 
For all $\theta \in (0, \infty )$ and all $d \ge 2$,
\begin{equation}\label{eq-pc}
	p_c^{(d)}(\theta)= m^{-1} + \dfrac{2d^2-1}{2(d-1)^2}\,m^{-2} + O(m^{-3}) +  O(m^{-1}V^{-1/3}),
\end{equation}
where the constants in the error terms may depend on $\theta$.
\end{thm}

\noindent
Observe that for $d \ge 4$, the correction term of order $m^{-2}$ is asymptotically larger than the 
width of the critical window, and that when $d =4,5,6$ the above expansion is again asymptotically 
precise, since we have $m^{-3} = O( m^{-1} V^{-1/3})$.  

To see the relevance of Theorem~\ref{thm-pc}, we compare it with other expansions of $p_c$ in the 
literature. The van der Hofstad and Slade \cite{HofSla06} proved that for percolation on $G$, with $G$ either the infinite 
lattice $\Z^d$ with nearest-neighbour edges or the hypercube $H(d,2)$, as $d \to \infty$, $p_c^G$ 
can be expanded up to three terms as
\begin{equation}
	p_c^{G} =  m^{-1} + m^{-2} + \frac{7}{2} m^{-3} + O(m^{-4}),
\end{equation}
where in both cases $m$ denotes the degree of the graph $G$. Moreover, they \cite{HofSla05}  also  
proved that, for any $N \in \bbN$,
\begin{equation}\label{e:latticeexp}
	p_c^{\Zd} = \sum_{k=1}^N a_k (2d)^{-k} + O((2d)^{-N-1}),
	\qquad p_c^{H(d,2)}(\theta) = \sum_{k=1}^N b_k d^{-k} + O(d^{-N-1}),
\end{equation}
where $(a_k),(b_k)$ are \emph{rational} coefficients. The critical window of the hypercube has width 
$O(d^{-1} 2^{-d/3})$, so we believe that the expansion cannot be asymptotically precise, regardless 
of the choice of $N$. Furthermore, it was conjectured that the expansion for $p_c^{\Zd}$, although it 
may exist, is \emph{divergent} for all $d$ as $N \to \infty$ (in the sense that the power series $z 
\mapsto \sum_{k=1}^\infty a_k z^k$ has radius of convergence $0$). 
\ch{We conjecture that the expansion for the Hamming graph is very different. We believe that for any $d \ge 2$ there exist coefficients $(c_k(d))$ such that}
\begin{equation}\label{e:pcconj}
	p_c^{\sss (d)}(\theta) = m^{-1} + \dfrac{2d^2-1}{2(d-1)^2}m^{-2} 
	+ \sum_{k=3}^{\lfloor d / 3 \rfloor} c_k(d) m^{-k} + O(m^{-1} V^{-1/3}),
\end{equation}
i.e., we conjecture that $p_c^{\sss (d)}$ has an asymptotically precise expansion in $m^{-1}$ of order 
$\lfloor d /3 \rfloor$ for all $d$. \ch{Heydenreich and van der Hofstad state the conjecture in \eqref{e:pcconj} 
as \cite[Open Problem~15.4]{HeyHof16}.}

Theorem~1.1 in \cite{BorChaHofSlaSpe05a} confirms this conjecture for $d=2,3$, and our current work confirms it for $d = 4,5,6$. The argument of van der Hofstad and Slade \cite{HofSla05} establishing \eqref{e:latticeexp} for the lattice and the hypercube crucially uses the fact that a ball of a radius $r$ restricted to a $d'$ dimensional subspace has the same shape for all $d\geq d'$, so that we can express each coefficient in terms of events that happen on a \emph{fixed} subgraph. Balls in the Hamming graph instead grow very rapidly when $n$ increases. Each coefficient is obtained as a limit and it will be more involved to prove the existence of this limit. Hence we do not have significant evidence suggesting that all coefficients in \eqref{e:pcconj} have to be rational.
\medskip

We note that the existence of a finite asymptotically precise expansion makes the proof of the critical 
window of the Hamming graph more challenging than for the hypercube. Roughly speaking, because 
the critical window of the hypercube is exponentially narrower than any of the expansion terms, we 
can approximate $p_c^{H(d,2)}$ up to any fixed order by a value $p$ that is in fact \emph{subcritical}, 
by choosing a negative coefficient for the error term. This allows \ch{one} to exploit the fact that $\chi (p)$ is 
poly-logarithmic in $V$, which simplifies the analysis considerably. In our case, the approximating $p$ 
will be much closer to $p_c^{\sss (d)}$, and so we need a much more refined analysis. We will explain 
this in more detail in Section~\ref{subsect:Pin}.


\subsection{Scaling limit of largest cluster sizes}

Besides offering an interesting comparison with other graphs with sufficiently weak geometry, the expansion of $p_c^{\sss (d)}$ also has another motivation. The Hamming graph is an excellent example to investigate the universality class of the ERRG, since it has a non-trivial geometry yet is highly mean field. See \cite{HofLuc10,HofLucSpe10,FedHofHul15,MilSen16} for a small sample of the literature from this perspective. A crucial motivation for the present paper is that it serves as a companion paper to \cite{FedHofHolHul16a}, where we establish the \emph{scaling limit} of the cluster sizes of the largest clusters within the critical window. More precisely, writing $\cluster_j$ for the $j$-th largest cluster, we prove that for any fixed $N \in \bbN$ and for $d=2,3,4$ the largest critical clusters of Hamming graph percolation satisfy
\begin{equation}\label{e:scaling}
	\big(V^{-2/3} |\cluster_j|\big) _{j \geq 1}
	\stackrel{\mathrm{d}}{\longrightarrow} (X_j)_{j\geq 1}
\end{equation}
for a certain sequence of \ch{$\theta$-dependent} continuous random variables $(X_i)_{i \in \bbN}$ supported on $[0,\infty)$. Aldous \cite{Ald97} proved this scaling limit for the ERRG. Since then, many other \ch{random graph} models have been shown to have the same (or at least a similar) scaling limit. \ch{See for instance \cite{BhaHofLee10a, BhaHofLee12,Jos14,NacPer10} and the references therein.} 
The above result for the Hamming graph, however, is the first indication that the same scaling occurs for models with an underlying high-dimensional geometry. Moreover, it is the most precise determination to date of the critical behaviour of percolation on a finite transitive graph (other than the ERRG scaling limit of Aldous). The proof of \eqref{e:scaling} and various other results in \cite{FedHofHolHul16a} crucially rely on the asymptotically precise determination of the critical window \ch{that we give here.} 


\ch{
\subsection{Alternative definition of the critical point}
It is worth noting that a disadvantage of the definition $p_c$ in \eqref{e:pcdef} is that it imposes an ad hoc relation between $p_c$ and $V^{1/3}$, which is known not to hold in general and believed to be associated with ``high-dimensional'' models. In other words, \eqref{e:pcdef} is possibly only a valid definition of $p_c$ for percolation models in the universality class of the ERRG. Nachmias and Peres in \cite{NacPer08} observed that it would be desirable to have a definition of $p_c$ that applies more generally, and they proposed
\begin{equation}
	\label{e:NacPerpc}
	\tilde p_c^G := \underset{p \in (0,1)}{\operatorname{argmax}} \frac{\frac{\mathrm{d}}{\mathrm{d} p} \chi_G(p)}{\chi_G(p)}
\end{equation}
as a definition of the critical point for \emph{any} graph $G$. Their motivation for this definition is that Russo's formula \cite{Rus81} implies that $p = \tilde p_c^G$ is the point where a small change in $p$ has the greatest impact on the relative size of the connected components, i.e., $\chi(p)$ changes most dramatically at $\tilde p_c^G$. 
A serious downside of this definition appears to be that $\tilde p_c^G$ may be very difficult to compute. Thus far, the only non-trivial determination of $\tilde p_c^G$ is given in recent work by Janson and Warnke \cite{JanWer16}. They determine that, for the ERRG, $|\tilde p_c^{K_n} - 1/n| = O(n^{-4/3})$, so $\tilde p_c^{K_n}$ is a point inside the critical window \eqref{e:critwindowERRG}, that $\chi(p)^{-1}\frac{\mathrm{d}}{\mathrm{d} p} \chi(p)$ around $\tilde p_c^{K_n}$ describes the critical window \eqref{e:critwindowERRG} as well, and that, interestingly, $\tilde p_c^{K_n}$ does \emph{not} equal either $1/n$ or $1/(n-1)$.
It would be interesting to see whether their methods can be applied to the current setting of percolation on $H(d,n)$.
}

\subsection{Susceptibility of the subcritical ERRG}

In Section \ref{sect-sus} we prove Theorem~\ref{thm-pc}, and also derive refined asymptotics for the susceptibility of a subcritical ERRG, \ch{its second moment, and its \emph{surplus:} given a connected graph $G$, let $\Spr(G) := |\Ecal (G)| -|\Vcal(G)| + 1$ denote the number of \emph{surplus edges} in $G$. Besides being interesting in their own right, these} will be crucial for proving the lower bound on $p_c^{\sss (d)}$, because the restriction of critical percolation on $H(d,n)$ to a one-dimensional subspace of $H(d,n)$ is equivalent to a subcritical ERRG. To prove the lower bound of Theorem~\ref{thm-pc} we rely on the following asymptotics, 
which, to the best of our knowledge, are sharper than results in the literature:

\begin{thm}[Second order asymptotics for susceptibility of the subcritical ERRG]
\label{secondorder} 
Let $G=G(n,p)$ be the ERRG with $p=\frac{\lambda}{n-1}$ and $0<\lambda<1$. Then as $n\to\infty$,
\begin{align}
\chi_G (p) \ch{= \E_p[|\cluster(v)|]} & = \dfrac{1}{1-\lambda}- \dfrac{2\lambda^2-\lambda^4}{2(1-\lambda)^4} n^{-1}  + O(n^{-2}),\label{e:susc}\\
\E_p [ |\cluster (v)|^2]& =  \dfrac{1}{(1-\lambda)^3} + O(n^{-1})
 ,\label{e:secondmoment}\\
 \ch{\E_p[\Spr(\cluster(v))] }&\ch{ = \frac{\lambda^3}{2(1-\lambda)^2} n^{-1} + O(n^{-2}).}\label{e:surplusbd}
\end{align}
\end{thm}

The second-order coefficient computed in \eqref{e:susc} improves the result by Durrett in \cite[Theorem 2.2.1]{Dur07},  which 
states that $\chi_G (p)=(1-\lambda)^{-1}-O(n^{-1})$, while \eqref{e:secondmoment} provides the 
matching lower bound to \ch{well-known upper bound derived with} the usual branching process domination. To achieve the sharper asymptotics 
we need a new way to encode the usual breadth-first search in the ERRG with the help of a branching 
random walk. We believe that there exists an infinite polynomial expansion of $\chi_G (p)$ in powers 
of $p$ for all $p=\frac{\lambda}{n-1}$ with $0<\lambda<1$. \ch{There is substantial literature related to \eqref{e:surplusbd}, see e.g.\
the classic book on random graphs by Bollob\'as \cite[Section 5.2]{Boll01} as well as the seminal paper by Janson, Knuth \L uczak and Pittel \cite{JanKnuLucPit93} 
computing generating functions of components having various cycle structures. As far as we are aware, the second order asymptotics in \eqref{e:surplusbd} is new.}


\subsection{Outline}

We prove Theorem~\ref{thm-pc} by separately proving a lower bound and an upper bound on 
$p_c^{\sss (d)}$. In Section~\ref{sect-sus} we prove Theorem~\ref{secondorder}. \ch{This theorem 
is used in Section~\ref{sect-low} to prove the lower bound in Theorem~\ref{thm-pc} with the help 
of an \emph{exploration process} that uses the fact that the restriction of critical bond percolation 
on $H(d,n)$ to a one-dimensional subspace has the same distribution as a subcritical ERRG. 
This is used to obtain a sharp enough \emph{branching process} upper bound on the susceptibility.}
In Section~\ref{sect-conn} we estimate connection probabilities and estimate bubble, triangle and 
polygon diagrams. In Section~\ref{sect-up} we prove the upper bound in Theorem~\ref{thm-pc} with the 
help of the \emph{lace expansion}. Perhaps surprisingly, these disparate methods yield compatible 
bounds, due to the fact that both methods are asymptotically sharp. The lace expansion method 
may be improved to prove Theorem~\ref{thm-pc}, but this would be more difficult than our current 
proof and  less interesting. We do not see how the \ch{exploration process} proof could be improved to also prove the 
upper bound in Theorem~\ref{thm-pc}.

\medskip


\section{Susceptibility of the subcritical Erd\H{o}s-R\'enyi Random Graph}\label{sect-sus}

In this section we prove Theorem~\ref{secondorder}. To give our estimate of the expected size of a 
subcritical cluster, we couple a breadth-first exploration process of the cluster to a process related 
to a \emph{Branching Random Walk} (BRW). The breadth-first exploration exploration process is 
defined in Section~\ref{ssect-BFe}, the branching random walk exploration in \ref{ssect-brwe}. 
The proof of the susceptibility asymptoticis is given in Section~\ref{ssect-sa}.


\subsection{Breadth-first/surplus exploration}
\label{ssect-BFe}
 
We start by defining a version of the breadth-first (BF) exploration. This is a very standard tool in the study 
of the ERRG (see e.g.\ \cite[Section~5.2.1]{Hofs17}). In a nutshell, a breadth-first exploration is a process 
that, starting from a vertex $v$, ``discovers'' its adjacent edges, ``activating'' the direct neighbours of $v$ 
in some fixed order, and then explores those vertices, discovering their adjacent edges and activating 
any unexplored, unactivated neighbours, and so on, always choosing the vertex that was activated the 
longest time ago as the next vertex to explore from. The BF exploration keeps track of which vertices 
have been explored (the ``dead'' set), which vertices have been activated but not explored (the ``active'' 
set), and the time at which a vertex was activated or explored. Crucially, the ``traditional'' BF exploration 
will only explore a vertex once, so the process terminates once all vertices are explored, and the edges 
associated with newly activated vertices describe a \emph{subtree} of the component of $v$, but the 
process provides little information about the \emph{surplus}, i.e., the discovered edges that do not activate 
new vertices (also sometimes referred to as the ``tree excess'' of the graph). For our purposes it is 
important that we also know about the surplus, so we consider the following modification of the BF 
exploration:

\begin{defi}[BF exploration process of a graph]
\label{def:BF} 
Given a graph $G = (\Vcal, \Ecal)$ and a vertex $v \in \Vcal$ we define the \emph{breadth-first/surplus 
(BF) exploration process} as the sequence of \emph{dead, active} and \emph{surplus sets} $(\Dbf (t), 
\Abf (t), \Sp (t))_{t \ge 0}$ as follows:
\begin{itemize}
	\item \textbf{Initiation.}
		Initiate the exploration with the dead, active and surplus sets at time $t=0$ as
		\begin{equation}
			\Dbf (0) := \varnothing, \qquad \Abf (0) := \{ v\}, \qquad  
			\Sp (0) := \varnothing,
		\end{equation}
		and at time $t=1$ as
		\begin{equation}
			\begin{aligned}
				\Dbf (1) &:= \{v\},\\ 
				\Abf (1) &:= \{w \, : \, \{v, w\} \in  \Ecal\},\\ 
				\Sp (1) &:= \varnothing.
			\end{aligned}
		\end{equation}
	\item \textbf{Time $t \geq 2$.}
	Choose the vertex $v^t \in \mathbf{A} (t-1)$ that minimizes $\min \{  i : v^t \in \mathbf{A} (i) \}$, 
	breaking ties according to an arbitrary but predetermined rule.\footnote{An example of such 
	a rule:  Fix an order on the vertex set $\Vcal$. If at step $t-1$ we have explored and/or activated 
	a total of $k$ vertices, and we activate $\ell$ more at step $t$, then we assign to these $\ell$ 
	newly explored vertices the labels $k+1$ through $k+\ell$, according to the order on $\Vcal$. 
	At time $s+1$ we explore from the active vertex with the smallest label.}
Update the active, dead and surplus sets as follows:
	\begin{equation}
		\begin{aligned}
			\Dbf (t) &:= \Dbf (t-1) \cup \{v^t\},\\ 
			\Abf (t) &:= (\mathbf{A} (t-1) \setminus \{v^t\}) \cup \{w \notin \Abf (t-1)
			\cup \Dbf (t-1) : \{v^{t}, w\} \in  \Ecal\},\\ 
			\Sp (t) &:= \Sp (t-1) \cup \{ \{v^t,w\} \in \Ecal : w\in \mathbf{A} (t-1)\}.
		\end{aligned}
	\end{equation}
	\item {\bf Stop.} Terminate the exploration when $\Abf (t) = \varnothing$. Set $T=t$.
\end{itemize}
\end{defi}

\noindent
Note that $\Dbf(t)$ and $\Abf(t)$ are subsets of $\Vcal$, whereas $\Sp(t)$ is a subset of $\Ecal$.
When $\Abf (t) = \varnothing$, this means that we have completely explored the connected component 
$\cluster (v)$ and $T = |\Dbf (T)| = |\cluster (v)|$. In the BF we find a new edge 
every time we activate a vertex (except the initial vertex $v$) or we discover an edge between active 
vertices. It follows that $|\Ecal (\cluster(v))|= |\Dbf (T)| -1 + |\Sp (T)|$. We conclude that that $|\Sp (T)|
=|\Ecal (\cluster(v))|- |\Dbf (T)|+1=\Spr (\cluster (v))$.


\subsection{The branching random walk exploration}
\label{ssect-brwe}

The subtree generated by a ``traditional'' BF exploration is often studied through a comparison to a 
branching process (see e.g.\ \cite{Hofs17,Dur07}). To study our BF exploration, we define a suitable 
extension, the \emph{branching random walk (BRW) exploration,} in which we randomly embed a 
branching process in the graph, and keep track of its self-intersections.\footnote{From now on the term 
\emph{nodes} will refer to elements of GW trees, while \emph{vertices} will refer to elements of graphs. 
Moreover, the \emph{progeny} of a node $x$ will indicate the set of vertices whose path to the root 
$\rho$ passes through $x$, while the \emph{children} of $x$ are only the vertices for which $x$ is the 
first vertex encountered on such a path. We write $\Cbf(x)$ for the set of children of $x$ in $\Tcal$.} 
This is made precise in the following definition:

\begin{defi}[Branching random walk]\label{def:brw}
Given an $m$-regular graph $G_m = (\Vcal, \Ecal)$ and $p \in [0,1]$, we define the 
$p$-\emph{branching random walk} ($p$-BRW) on $G_m$ started at $v \in \Vcal$ as 
the pair $(\Tcal, \phi_v)$, where $\Tcal$ is a $\Bin(m,p)$ Galton-Watson tree, and $\phi_v$ 
is a random mapping of $\Tcal$ into the vertex set $\Vcal$ whose law satisfies: (1) $\phi_v$ 
maps the root $\rho$ of $ \Tcal$ to $v$; (2) given any node $x \in \Tcal$ and its set of children 
$\Cbf (x) \subset \Tcal$, the marginal law of $\phi_v(\Cbf (x))$ is the same as that of $|\Cbf(x)|$ 
\emph{distinct} neighbours of $\phi_{v} (x)$ in $G_m$ chosen uniformly at random, independently 
for all $x \in \Tcal$. (Here, for a set $A \subset \Tcal$ and a mapping $\phi_v\colon \Tcal\to G$, we 
define $\phi_v(A) = \cup_{a \in A} \phi_v(a)$, and by convention set $\phi_v(\varnothing) 
= \varnothing$.)
\end{defi}

Next, we define a process that explores a $p$-BRW and keeps track of any self-intersections. Briefly, 
the idea is that we explore the $p$-BRW by exploring the tree $\Tcal$ in a breadth-first fashion from 
the root upward. If the $p$-BRW intersects its own trace, then we declare the particle that intersected, 
and all its offspring, to have become ``ghosts''. We differentiate between particles that became ghosts 
through intersecting with active and dead vertices. In Proposition~\ref{coupling} below we prove that 
this exploration process can be coupled to a BF exploration of a percolation cluster:

\begin{defi}[BRW exploration process]
\label{def:brwe} 
Given an $m$-regular graph $G_m = (\Vcal, \Ecal)$, a vertex $v \in \Vcal$, and a $p$-BRW 
$(\Tcal, \phi_v)$ on $G_m$, we define the \emph{BRW exploration process} $( \Att  (t), \Dtt  (t),
\Ptt ^{\Att} (t), \Ptt ^{\Dtt}(t))_{t = 0}^T$ as the sequence of \emph{dead, active, active ghost} 
and \emph{dead ghost} sets as follows:
\begin{itemize}
	\item \textbf{Initiation.}
		Initiate the exploration with the \emph{dead, active, active ghost} and 
		\emph{dead ghost} sets at time $t=0$ as
		\begin{equation}
			\Dtt  (0) := \varnothing, \quad \Att  (0) := \{ \rho\}, \quad 
			\Ptt ^{\Att} (0)=\varnothing , \quad \Ptt ^\Dtt(0)	=\varnothing,
		\end{equation}
		and at time $t=1$ as
		\begin{equation}
			\Dtt  (1) := \{\rho\}, \quad
			\Att  (1) := \{y\in \Cbf (\rho )\}, \quad
			\Ptt ^{\Att} (1) := \varnothing , \quad
			\Ptt ^{\Dtt}(1)	:= \varnothing .
		\end{equation}

	\item \textbf{Time $t \geq 2$.} 
		Choose the node  $x^t \in \Att  (t-1)$ that minimizes $\min \{  i : x^t \in \Att (i) \}$, 
		breaking ties according to an arbitrary but predetermined rule, and update the 
		exploration as follows:
		\begin{equation}\label{e:BRWdeftbig}
			\begin{aligned}
				\Dtt  (t) &:= \Dtt  (t-1) \cup \{x^t\},\\ 
				\Att  (t) &:= (\Att  (t-1) \setminus \{x^t\}) \cup \big\{y\in \Cbf (x^t): \phi_{v}(y)
				\notin \phi_{v}\big(\Dtt  (t-1) \cup \Att   (t-1)\big) \big\},\\
				\Ptt ^{\Att} (t) & :=  \Ptt ^{\Att}(t-1) \cup \big\{ y\in \Cbf(x^t): \phi_{v}(y) 
				\in \phi_{v}\big(\Att (t-1)\big)\big\},\\
				\Ptt ^{\Dtt}(t) & := \Ptt ^{\Dtt}(t-1) \cup \big\{ y \in \Cbf (x^t): \phi_{v}(y) 
				\in \phi_{v}\big(\Dtt (t-1)\big)\big\}.
			\end{aligned}
		\end{equation}
	\item \textbf{Stop.}
		If $\Att  (t)= \varnothing$, then terminate the exploration. Set $T=t$.
\end{itemize}
\end{defi}

Using the BRW exploration, we define the subgraph $\tilde \cluster(v)$ as the graph traced out by a 
$p$-BRW where the particles are killed when they intersect with the active set. More precisely, we 
let $\tilde \Tcal$ be the subtree in $\Tcal$ induced by $\Dtt(T) \cup \Ptt^{\Att}(T)$, and define 
\begin{equation}\label{e:Ctildedef}
	\tilde \cluster(v) := \big( \phi_v(\Dtt(T)), \{\{\phi_v(x), \phi_v(y)\} : \{x,y\} \in \tilde \Tcal\} \big).
\end{equation}
Note that, by Definition~\ref{def:brwe}, $\phi_v(\Dtt(T) \cup \Ptt^{\Att}(T)) = \phi_v(\Dtt(T))$, so 
$\tilde \cluster(v)$ is indeed a subgraph of $G_m = (\Vcal, \Ecal)$. 

We now show that $\tilde \cluster(v)$ has the same law as $\cluster(v)$, the connected component of 
$v$ in an ERRG, by coupling the BF and BRW explorations:

\begin{prop}[Coupling of BF and BRW explorations]
\label{coupling} 
Consider percolation on an $m$-regular graph $G_m$ with parameter $p$. Consider the BF exploration 
on the percolated graph $G_m(p)$ and the $p$-BRW exploration processes on $G_m$, both starting 
from the vertex $v$ (and using the same tie-breaking rule). Then $\cluster (v)$ with respect to $\P_p$ 
has the same law as $\tilde{\cluster}(v)$.
\end{prop}

\proof
We show inductively that we can couple each step of the BRW and of the BF exploration in such a 
way that $\tilde\cluster(v)=\cluster (v)$ almost surely. We start by showing that there exists a coupling such 
that for all $t \ge 0$,
\begin{equation}\label{couple}
	\begin{aligned}
		\Dbf (t)&= \phi_v (\Dtt  (t)), \\
		\Abf (t)&= \phi_v (\Att  (t)), \\
		\Sp (t)&= \bigcup_{s \leq t} \big\{ \{\phi_v(x^{s}),w\}:
		w \in \phi_v(\Ptt ^{\Att} (s)\setminus \Ptt ^{\Att} (s-1))\big\}.
	\end{aligned}
\end{equation}

We start with the inductive base. At time $t=0$, by Definitions \ref{def:BF} and \ref{def:brwe},
\begin{equation}
	\begin{aligned}
		\Dbf (0)&= \varnothing = \phi_v (\Dtt  (0)),\\
		\Abf (0) &= \{v\} = \phi_v (\{\rho\}) =  \phi_v (\Att  (0)),\\
		\Sp (0)&=\varnothing = \big\{ \{\phi_v(x^{0}),w\}: w \in \phi_v(\Ptt ^{\Att} (0))\big\} .
	\end{aligned}
\end{equation}
Next, we prove the inductive step: the induction hypothesis is that the relations in \eqref{couple} 
holds for all $r<t$. We extend the coupling so that it also holds at time $t$. Our assumption is that 
we use the same tie-breaking rule for both explorations, so by the induction hypothesis we choose 
$v^t=\phi _{v} (x^t)$.

Given $\phi_{v}(x^t)$, fix a set $\mathcal U _k^t =\{ u^1, u^2,\ldots,u^k \}$ of $k$ neighbours 
of $\phi_{v}(x^t)$. By Definition~\ref{def:brw}, the mapping $\phi_{v}$ is such that $|\Cbf (x^t)|$ 
neighbours of $\phi_{v}(x^t)$ are distinct neighbours chosen uniformly at random, so
\begin{equation}\label{eqX}
	\begin{aligned}
		\P \big(\phi_{v} (\Cbf (x^t)) = \mathcal U _k^t \big) &= \P ( |\Cbf (x^t)| =k)\,
		\P\big(\phi_{v} (\Cbf (x^t)) = \mathcal U _k^t \mid |\Cbf (x^t)| =k\big)\\ 
		&= {{m}\choose{k}} p^k (1-p)^{m-k} {{m}\choose{k}} ^{-1} = p^k (1-p)^{m-k}.
	\end{aligned}
\end{equation}

Next, consider the BF exploration at time $t$. Given $(\Dbf(s), \Abf(s), \Sp(s))_{s=0}^{t-1}$, we 
can determine $v^t$. For $t \ge 0$, let $\Nbf(t)$ denote an independent set-valued random variable 
that contains the vertex $w$ with probability $p$, independently for all $w$ such that $\{v^t,w\} \in 
\Ecal$, so that $\P \big(\Nbf (t)=\mathcal U _k^t \big) = p^k (1-p)^{m-k}$. For every set $\Ucal_k^t$ 
of neighbours of $v^t$ we have $\P \big(\Nbf (t)=\mathcal U _k^t \big)=\P \big(\phi_{v} (\Cbf (x^t)) 
= \mathcal U _k^t \big)$, and so there exists a trivial coupling of $\Nbf(t)$ and $\phi_v(\Cbf(x^t))$ 
such that $\P \big(\Nbf (t)=\phi_{v} (\Cbf (x^t)\big)=1$.

Consider an edge $\{v^t, w\}$. Observe that if $w \notin \Dbf(t-1)$, then $\{v^t,w\}$ has not been 
discovered by the exploration, so it is open in the percolation conditionally independently with probability~$p$, while if $w \in \Dbf(t-1)$, then $\{v^t, w\}$ has been discovered in the BF exploration, so its 
status can be determined from $(\Dbf(s), \Abf(s), \Sp(s))_{s=0}^{t-1}$. Let $\Xbf(t)$ denote the vertices 
that are end-points of edges that are discovered in the $t$-th step, i.e.,
\begin{equation}
	\Xbf(t) := \big\{\{v^t, w\} \, : \, \{v^t, w\} \in \cluster(v) \setminus
	 \{\{v^t, u\} \, : \, u \in \Dbf(t-1)\} \big\}.
\end{equation}
Note that if $w \in \Xbf(t)$, then either $w$ becomes activated at time $t$ or $w \in \Abf(t-1)$. By the 
above observation, we can couple $\Xbf(t)$ to $\Nbf(t)$ such that $\Xbf(t) = \Nbf(t) \setminus \Dbf(t-1)$ 
almost surely, conditionally on $(\Dbf(s), \Abf(s), \Sp(s))_{s=0}^{t-1}$. 

Consider henceforth the setting in which $\Xbf(t)$, $\Nbf(t)$, $\phi_v(\Cbf(v^t))$, and $(\Dbf(s), \Abf(s), 
\Sp(s))_{s=0}^{t-1}$ are \emph{simultaneously} coupled according to the above description. (Since 
both $\Nbf(t)$ and $\phi_v(\Cbf(x^t))$ are essentially independent $p$-random subsets, it is easy to 
make this coupling explicit; we leave those details to the reader.) Using this coupling, the induction 
hypothesis \eqref{couple}, and Definitions \ref{def:BF} and \ref{def:brwe}, we derive
\begin{equation}
	\Dbf (t)\setminus \Dbf (t-1) = \{v^t\}=\phi (x^t) = \phi_v (\Dtt  (t))\setminus \phi_v (\Dtt  (t-1)),  
\end{equation}
and
\begin{equation}
	\begin{aligned}
		\Abf (t)\setminus \Abf (t-1)
		& = \Xbf(t) \setminus \Abf(t-1)\\
		&= \Nbf (t) \setminus (\Dbf (t-1) \cup \Abf (t-1))\\
		&=(\phi_{v} (\Cbf (x^t)) \setminus (\phi_v(\Dtt (t-1)) \cup \phi_v(\Att (t-1)))\\
		&= \phi_v(\Att (t)) \setminus \phi_v(\Att (t-1)),
	\end{aligned}
\end{equation}
and
\begin{equation}
	\begin{aligned}
		\Sp (t)\setminus \Sp(t-1) &= \{\{v^t, w\} \, : \, w \in \Xbf(t) \cap \Abf(t-1)\}\\
		& = \{ \{v^t , w\} \, : \, w \in \Nbf (t) \cap \Abf (t-1)\}\\
		&= \{ \{\phi_v (x^t),w\} \, : \, w \in \phi_{v} (\Cbf (x^t)) \cap \phi_v(\Att (t-1))\}\\ 
		&= \{ \{\phi_v (x^t),w\} \, : \, w \in  \phi_v(\Ptt ^{\Att} (t)\setminus \Ptt ^{\Att} (t-1))\}.
	\end{aligned}
\end{equation}
Since $\phi_v(x^t)=v^t$, we obtain that \eqref{couple} holds also at time $t$ almost surely, and thus, 
by induction, for all $t \in \{0,1,\dots,T\}$, almost surely.

To conclude the proof, we show that the coupling \eqref{couple} for all $t \in \{0,1,\dots,T\}$ implies 
that $\tilde \cluster(v) = \cluster(v)$ almost surely. Recall the definition of $\tilde \cluster(v)$ in \eqref{e:Ctildedef}, 
and of $\tilde \Tcal$ above it. Since $\phi_v(\Dtt(T)) = \Dbf(T)$, it follows directly from \eqref{couple} 
that the vertex sets of $\tilde \cluster(v)$ and $\cluster(v)$ coincide. To see that the edge sets coincide, 
note that, by Definition~\ref{def:brwe}, $\tilde \Tcal$ contains only edges $\{x,y\} \in \Tcal$ such that 
$y\in \Cbf (x)$, with $x \in \Dtt(T)$ and $y \in \Dtt(T) \cup \Ptt^{\Att}(T)$. Indeed, by the construction 
of the BRW exploration it is impossible that both $x, y \in \Ptt^{\Att}(T)$, since vertices in $\Ptt^{\Att}$ 
are never explored further. Let $s = \min \{t: x \in \Dtt (t)\}$. Then, from the definition of the BRW 
exploration, $x=x^s$ and $y \in (\Att (s)\setminus \Att (s-1)) \cup (\Ptt^{\Att}(s)\setminus \Ptt^{\Att}(s-1))$. 
We then obtain
\begin{equation}
\Ecal (\tilde\cluster (v))= \bigcup_{s = 0}^{T} \big\{ \{\phi_v(x^s), \phi_v(y)\} 
\, : \, \phi_v(y) \in \phi_v((\Att (s)\setminus \Att (s-1))\cup(\Ptt^{\Att}(s) 
\setminus \Ptt^{\Att}(s-1)) \big\} .
\end{equation}
An application of \eqref{couple} now completes the proof. \qed
\medskip

We conclude by deriving some consequences of Proposition~\ref{coupling} that will be 
useful in the proof of Theorem~\ref{secondorder}:

\begin{cor}\label{cor-coup}
Consider a BF exploration on an $m$-regular graph $G_m$ with parameter $p$, and 
a $p$-BRW exploration processes on $G_m$, both starting from the vertex $v$. Then
\begin{align}
	| \Dtt (T) | &\overset{d}{=} |\cluster (v)|, \label{eqsize}\\
	|\Ptt ^{\Att} (T) | &\overset{d}{=} \Spr(\cluster (v)),\label{eqsurp}\\
	\E[|\Ptt ^{\Dtt}(T)|]&= p \E \Big[ \sum_{t=1}^T |\Dtt  (t-1)|\Big].\label{eqdead}
\end{align}
\end{cor}

\proof
Equations \eqref{eqsize} and \eqref{eqsurp} follow immediately from Proposition~\ref{coupling}. 
To prove \eqref{eqdead} we use \ch{\eqref{e:BRWdeftbig} and} \eqref{eqX}, from which we get 
\begin{equation}
	\E [|\Ptt ^{\Dtt}(t)\setminus \Ptt ^{\Dtt}(t-1)|\mid \Fcal (t-1)] =p|\Dtt (t-1)|,
\end{equation}
where $\Fcal (t-1)$ is the $\sigma$-algebra generated by the first $t-1$ steps of the BRW exploration. 
Summing over $t \leq T$ we get the claim. \qed


\subsection{Proof of the susceptibility asymptotics}
\label{ssect-sa}

\proof[Proof of Theorem \ref{secondorder}]
Consider a BF exploration on $G(n,p)$ and a BRW exploration of a 
$p$-BRW on $K_n$, \ch{both with $p = \frac{\lambda}{n-1}$,} coupled as in the proof of Proposition~\ref{coupling}. Corollary~\ref{cor-coup} 
implies that $\chi_G (p)= \E [|\Cctv|]=\E [|\Dtt (T)|]$. We know that $\E [|\Tcal | ] = \frac{1}{1-\lambda}$ 
(see \cite{Pak71}). From the definition of the BRW exploration we know that $\Tcal \setminus \Dtt (T)$ 
consists of the nodes $\Ptt ^{\Dtt}(T) \cup \Ptt ^{\Att}(T)$ and their progeny in $\Tcal$. Using \eqref{eqsize} 
and \eqref{eqsurp}, we can thus write
\begin{equation}
	\begin{aligned}
		\E_p[|\cluster (v)|]&=  \E[|\Tcal |]-\E[|\Tcal  |]\,\E[|\Ptt ^{\Att}(T) \cup \Ptt ^{\Dtt}(T)|]\\
		&= \E[|\Tcal |]-\E[|\Tcal  |]\,\Big(\E_p[\Spr(\cluster (v))]
		+p\E\Big[ \sum_{t=1}^T |\Dtt  (t-1)|\Big]\Big).
	\end{aligned}
\end{equation}
Since $|\Dtt  (t-1)|=t-1$ for all $t \leq T$ and $|\Dtt(T)|\overset{d}{=}|\cluster(v)|$,
\begin{equation}
	\E\Big[\sum_{t=1}^T |\Dtt  (t-1)|\Big]=\tfrac{1}{2}\E[ T(T-1)]
	=\tfrac{1}{2}\E \big[ |\cluster (v)|^2-|\cluster (v)|\big].
\end{equation}
As a result we obtain that
\begin{equation}\label{e:Cibd}
	\E_p[|\cluster (v)|]= \E[|\Tcal |] \left(1-\Ep[\Spr(\cluster (v))]
	-\tfrac12 p \E_p\big[|\cluster (v)|^2-|\cluster (v)|\big]\right).
\end{equation}
Pakes proved in \cite[Section 2.2]{Pak71} that 
\begin{equation}
	\Var (|\Tcal |)
	= \frac{\lambda}{(1-\lambda)^3}+O(p),
\end{equation}
so that $\E [|\Tcal |^2]= \ch{(1-\lambda)^{-3}}+O(p)$. Moreover, Durrett shows in \cite[Section 2]{Dur07} that $\E_p[|\cluster (v)|]=\ch{(1-\lambda)^{-1}}-O(n^{-1})$.

We will prove that $\E_p[|\cluster (v)|^2]= \ch{(1-\lambda)^{-3}}+O(p)$. It follows from 
Proposition~\ref{coupling} that $ \E_p[|\cluster (v)|^2]= \E[|\Dtt (T)|^2] \leq \E [|\Tcal |^2]$, which 
establishes the upper bound. To determine the lower bound, we write
\begin{equation}\label{e:Ccalvsquared}
	\E_p[|\cluster (v)|^2] = \E [ (|\Tcal |-|\Tcal \setminus \Dtt(T) |)^2] 
	\geq \E [|\Tcal|^2] - 2 \E [ |\Tcal||\Tcal\setminus \Dtt(T) |],
\end{equation}
so it remains to prove that $\E [ |\Tcal ||\Tcal \setminus \Dtt (T) |]= O(p)$. We write
\begin{equation}\label{k4}
	\begin{aligned}
		\E [ |\Tcal ||\Tcal \setminus \Dtt (T)|]  
		&=\E \Big[ \sum_{x,y \in \Tcal } \indi_{\{x \notin \Dtt (T) \}}\Big]\\ 
		& = \sum_{k=1}^{\infty} \P (|\Tcal| =k)\,k\,\E \Big[ \sum_{x } \indi_{\{x \notin \Dtt (T), 
		x \in \Tcal \}} \big| |\Tcal |=k\Big].
	\end{aligned}
\end{equation}
Now we bound $\E \big[ \indi_{\{x \notin \Dtt (T), x \in \Tcal \}} \big| |\Tcal |=k\big]$. Suppose that 
$d_\Tcal(\rho,x) = L$, and that $\eta$ is the path in $\Tcal$ with $\eta(0) = \rho$ and $\eta(L) =x$.
From the description of the BRW exploration it follows that $x \notin \Dtt (T) $ if and only if there 
exists a $t \in \{0, \dots, T\}$ such that $v^t \in \eta$ and $\phi_v(v^t) \in \Ptt^\Att (t-1) \cup \Ptt^\Dtt (t-1)$. 
Since the mapping of children of a node in the BRW is done uniformly at random, $\phi_v(\eta)$ has 
the same distribution as a simple random walk path on $K_n$, conditioned on going from $\rho$ to 
$x$ in $L$ steps. Conditionally on $|\Tcal|=k$, we have $|\Att (t-1)\cup \Dtt (t-1)| \leq k$ for all $t$, 
so for all $t$ the probability that $\phi_v(v^t) \in \phi_v(\Att (t-1)\cup \Dtt (t-1))$ is at most $k/(n-1)$ 
by the symmetry of $K_n$. Conditionally on $|\Tcal|=k$, we have $L \le k-1$, and so
\begin{equation}
	\E \big[ \indi_{\{x \notin \Dtt (T), x \in \Tcal\}} \big| |\Tcal|=k\big] \leq \frac{k (k-1)}{n-1}.
\end{equation}
Inserting this into \eqref{k4} we obtain
\begin{equation}
	\E [ |\Tcal ||\Tcal \setminus \Dtt (T) |] 
	= \sum_{k=1}^{n} \P (|\Tcal| =k)\,k^2\,\frac{k(k-1)}{n-1}\leq \dfrac{p}{\lambda} \E [|\Tcal|^4] 
	= O(n^{-1}).
\end{equation}
Inserting this into \eqref{e:Ccalvsquared} we get \eqref{e:secondmoment}, and thus it also
follows that
\begin{equation}\label{e:Csecmombd}
	\E_p[|\cluster (v)|^2]-\E_p[|\cluster (v)|]
	= \frac{1}{(1-\lambda)^3} - \frac{1}{1-\lambda} +O(n^{-1})
	=\dfrac{2\lambda - \lambda^2}{(1-\lambda)^3}+O(n^{-1}).  
\end{equation}

Next, we compute $\E_p[\Spr(\cluster(v))]$. Note that $\Spr(\cluster(v))$ is bounded from above by the 
number of vertex-disjoint cycles in $\cluster(v)$, since removing a surplus edge from a graph destroys 
at least one such cycle. Here, given a graph $G = (\Vcal, \Ecal)$ and $v_1,\dots, v_k \in \Vcal$, we say that a 
subgraph $L_k$ is a \ch{(vertex-disjoint)} \emph{cycle of length $k$} if $L_k$ has vertex set $\{ v^1, \dots, v^k\} \subset 
\Vcal$ and edge set $\cup_{i=1}^k \{v^i, v^{i + 1 \,\mathrm{mod}\, k}\}$. We write $\Lcal_k(G)$ for 
the set of cycles of length $k$ in a graph $G$. For $G(n,p)$ we may thus bound
\begin{equation}
	\E_p [\Spr(\cluster(v))]\leq \sum_{k=3}^n \sum_{L_k \in \Lcal_k(K_n)} 
	\P_p \big(L_k \in \Lcal_k(\cluster(v))\big).
\end{equation}
Note that $|\Lcal_k(K_n)| =\frac{n!}{2k(n-k)!}$. The probability that a given set of $k$ edges is open 
in $G(n,p)$ is $p^k$. Moreover, the probability that a vertex $w$ in $G(n,p)$ is connected to a given 
set of $k$ vertices is at most $ k \frac{\E_p [|\cluster (w)|]}{n} \leq \frac{k}{(1-\lambda)n}$ for any $w$, 
by the symmetry of $K_n$. Combining these estimates, and using that $p = \frac{\lambda}{n-1}$, we bound
\begin{equation}\label{e:Spbd}
	\begin{split}
		\E_p[\Spr(\cluster(v))]&\leq\sum_{k=3}^n \dfrac{n!}{2k(n-k)!}p^k 
		\dfrac{k }{(1-\lambda)n}(1+O(n^{-1}))\\
		&= \sum_{k=3}^n \frac{(n-1)!}{(n-k)!}\, \frac{1}{(n-1)^{k-1}}\, 
		\frac{\lambda^{k-1}}{2(1-\lambda)}\,p\,(1+O(n^{-1}))\\
		&\ch{\le \sum_{k=3}^n \dfrac{\lambda^{k}}{2(1-\lambda)} n^{-1}+O(n^{-2})
		= \dfrac{\lambda^3}{2(1-\lambda)^2} n^{-1}+O(n^{-2}).}
	\end{split}
\end{equation}
\ch{This establishes the upper bound in \eqref{e:surplusbd}.}

It remains to prove a matching lower bound. Before we start, let us recall a standard tool from percolation
theory: the \emph{van den Berg-Kesten (BK) inequality} \cite{BerKes85}: We say that an event $A$ is 
\emph{increasing} with respect to $p$ if $\P_p(A) \ge \P_q(A)$ whenever $p \ge q$. We say that two 
increasing events $A$ and $B$ \emph{occur disjointly}, and write $A \circ B$, if the occurrence of $A$ 
and $B$ can be verified by inspecting disjoint sets of edges (which may depend on the percolation
configuration). For instance, the event $\{v \leftrightarrow w\}$ is increasing, and $\{v \leftrightarrow w\} 
\circ \{v' \leftrightarrow w'\}$ implies that there exists a path of open edges between $v$ and $w$ and 
another path of open edges between $v'$ and $w'$, and that these paths are edge disjoint. The 
BK-inequality states that $\P_p(A \circ B) \le \P_p(A) \P_p(B)$. See e.g.\ Grimmett's classic book on percolation 
\cite{Grim99} for more details. 
\medskip

We use the BK-inequality to prove a lower bound on the expected surplus. Since the removal of a surplus 
edge must destroy at least one cycle in the graph, we can bound $\Spr(\cluster(v))$ from below by the 
number of \ch{vertex-disjoint cycles in $\cluster(v)$ that are edge-disjoint from any other cycle in $\cluster(v)$.} A cycle is edge-disjoint from other cycles if and only if the cycle does not contain a pair of vertices that are connected by a path outside the cycle. Each cycle of length $k$ has $k(k-1)/2$ pairs of vertices, and the probability that any two vertices are connected is at most $\frac{1}{(1-\lambda)n}$ by the symmetry of $K_n$. Writing $\Vcal(L_k)$ for the vertex set of the cycle $L_k$ and using the same reasoning as in \eqref{e:Spbd} as well as inclusion-exclusion, we thus 
obtain
\begin{align}\label{e:Splb}
		\E_p[\Spr(\cluster(v))]  & \ge \sum_{k=3}^{n} \sum_{L_k \in \Lcal_k(K_n)} 
		\bigg[\P_p \big(L_k \in \Lcal_k(\cluster(v))\big)\\
		& \qquad\qquad\qquad\qquad -  \P_p \Big(\bigcup_{\{x,y\} \subset \Vcal(L_k)}
		\big\{ L_k \in \Lcal_k(\cluster(v)) \big\} \circ  \{ x \leftrightarrow y \} \Big)\bigg]\nn\\
		& \ge \sum_{k=3}^n \frac{n!}{(n-k)!}\, \frac{p^{k}}{2(1-\lambda)n}\,  
		\Big(1 - \frac{k(k-1)}{2(1-\lambda)n}\Big)\nn\\
		&\ch{ \ge \sum_{k=3}^{n}\Big( \dfrac{\lambda^{k}}{2(1-\lambda)}\,n^{-1}
		+O(n^{-2})\Big)\Big(1- \dfrac{k(k-1)}{2 (1-\lambda)} n^{-1} \Big)}\nn\\
		&\ch{ = \dfrac{\lambda^3}{2(1-\lambda)^2}n^{-1} +O(n^{-2}),}\nn
	\end{align}
where in the first step the union is over all two-element subsets of $\Vcal(L_k)$, in the second step we 
first use the union bound and then the BK-inequality, and in the third step we use that $p = \lambda/(n-1)$ 
and $\frac{(n-1)!}{(n-k)!} = (n-1)^{k-1} - O(k n^{k-2})$. \ch{This completes the proof of \eqref{e:surplusbd}.
\medskip}

Inserting the bounds \eqref{e:Csecmombd}, \eqref{e:Spbd}, and \eqref{e:Splb} into \eqref{e:Cibd}, we 
conclude that
\begin{equation}
\begin{split}
	\E_p[|\cluster (v)|]
	&=\ch{\dfrac{1}{1-\lambda}\Big(1 -\frac{\lambda^3}{2(1-\lambda)^2} n^{-1} - \tfrac12 p \frac{2 \lambda -\lambda^2}{(1-\lambda)^3} + O(n^{-2}) \Big)}\\
	&= \dfrac{1}{1-\lambda}- \dfrac{2\lambda^2-\lambda^4}{2(1-\lambda)^4} n^{-1} + O(n^{-2}),
\end{split}
\end{equation}
\ch{which proves \eqref{e:susc} and thus completes the proof of Theorem~\ref{secondorder}.} 
\qed


\section{The lower bound on $p_c^{\sss (d)}(\theta)$ via an exploration process}\label{sect-low}

In this section we use the bound on the susceptibility of the subcritical ERRG to determine a 
lower bound on $p_c^{\sss (d)}(\theta)$, the critical value of the Hamming graph. We achieve 
this by bounding $\chi(p)$ from above with the use of an exploration process, and then substituting 
this bound into $\chi(p_c^{\sss (d)}) = \theta V^{-1/3}$, the equation that defines $p_c^{\sss (d)}$ 
(recall \eqref{e:pcdef}). The exploration process that we use is designed with the geometry of 
the Hamming graph in mind, so let us start by investigating this geometry further.

Recall that the Hamming graph $H(d,n)$ can be viewed as the $(d-1)$-fold Cartesian product 
of complete graphs $K_n$. If we arrange the $n^d$ vertices of $H(d,n)$ on a $d$-dimensional 
hypercubic grid in the obvious way, then the edges of $H(d,n)$ are precisely those edges that 
have both end-points on a line that is parallel to an axis of the grid. This inspires the following 
definition: Given some $i \in \{1,\dots, d\}$ and a vertex $v = (v_1, \dots, v_d ) \in \{0,\dots,n-1\}^d$, 
we call the subgraph of $H(d,n)$ induced by the set 
\begin{equation}
	\big\{(v_1,\dots, v_{i-1}, w, v_{i+1}, \dots, v_{d}) \in \Vcal \, : \, w \in \{0,\dots,n-1\} \big\},
\end{equation}
the $i$\emph{-directional line of $H(d,n)$ through $v$}. When $i$ and $v$ are unimportant we refer 
to such subgraphs simply as \emph{lines}. We write $\cluster_i(v)$ for the set of vertices that can be 
reached by a path of open edges in the $i$-directional line through $v$. Note that any line of $H(d,n)$ 
is isomorphic to $K_n$, so that $\cluster_i(v)$ has the same law as an Erd\H{o}s-R\'enyi random graph 
on $n$ vertices with parameter $p$. Moreover, because this is a graph on $n$ vertices and we choose 
$p =\frac{1+O(m^{-1})}{d(n-1)}$ with $d \ge 2$ (in accordance with \eqref{e:critwindowBCHHS}), this 
ERRG is \emph{subcritical.} Writing $\chil (p) : = \Ep[|\cluster_i(v)|]$ for the expected size of a connected
 component \emph{within a line} \ch{(i.e., the set over vertices that can be reached from $v$ using only open edges in the line),} we get from Theorem \ref{secondorder} that
\begin{equation}\label{e:chilpsimple}
	\chil  \Big(\frac{1+O(m^{-1})}{d(n-1)}\Big) = \frac{1}{1-p(d-1)} 
	\Big(1 - \dfrac{2d^2-1}{2(d-1)^3}m^{-1} + O(m^{-2})\Big).
\end{equation}
We use this fact repeatedly below.

We next define the exploration process that allows us to estimate $\chi(p)$. To first order, this estimation 
simply yields a Galton-Watson branching process. But this is an overestimate, and we can give a 
(negative) second order correction to it by correcting for the over counting that arises because we 
ignored loops in the graph. We thus will have to find a bound on the number of loops. An important 
insight into the structure of percolation on the Hamming graph is that loops are much more likely to 
occur within lines than outside lines. Our exploration crucially uses this fact: we only subtract the 
correction for loops within lines, which gives us the desired upper bound.

Roughly speaking, the line-wise exploration process defined below works as follows. We have two 
sets, the \emph{active} and \emph{dead} sets $\Asf$ and $\Dsf$. We start with a single vertex in the 
active set. At any given time we move an active vertex to the dead set, and add all the vertices connected 
to that vertex \emph{through a line} to the active set. Because we want to avoid ``feedback loops'' in 
the process, we need to keep track of the line that we have previously explored from. The \emph{parent 
set} $\Psf$ of ordered pairs of vertices and their parents in the exploration is a technical addition to 
the process that takes care of this. The process stops when the active set becomes empty. 

\begin{defi}
\label{def-linewise}
The \emph{line-wise exploration process} $(\Asf (t) , \Dsf (t), \Psf(t))_{t= 0}^T$ on $H(d,n)$ started at 
the vertex $v$ is the $T$-step discrete-time process defined as follows:
\begin{itemize}
	\item \textbf{Initiation.}  Define the \emph{dead, active} and \emph{parent sets} 
	at time $t=0$ 
	\begin{equation*}
		\Dsf (0) := \varnothing, \qquad \Asf (0) := \{ v\}, \qquad \Psf(0) := \varnothing. 
	\end{equation*}
	and at time $t=1$
	\[
		\Dsf (1) : = \{v\},  \qquad \Asf (1):= \bigcup_{i\in [d]} \{w : w \in \cluster_i (v)\} 
		\setminus \{v\}, \qquad \Psf(1) :=  \big\{\{ w, v\} \in \Ecal : w \in \Asf(1)\}.
	\]
	\item  \textbf{Time $t \geq 2$.} Choose a vertex $v^t$ according to an arbitrary but 
	predetermined rule from $\Asf(t-1)$. Let $u$ be the vertex such that 
	$\{ u, v^t\} \in \Psf(t-1)$ and write $j_t$ for the unique direction such that 
	$v^t_{j_t} \neq u_{j_t}$. Update
	\[
	\begin{split}
		\Dsf (t) &:=\Dsf (t-1) \cup  \{v^t \}, \\ \Asf (t)&:=\Big( \Asf (t-1) \cup 
		\bigcup_{i\in [d]\colon i \neq j_t} \big\{w : w \in \cluster_i (v^t)\big\} \Big)
		 \setminus \{v^t\},\\
		\Psf(t) &:= \Psf(t-1) \cup \big\{\{w,v^t\} \in \Ecal : w \in  \Asf(t) \setminus 
		\Asf(t-1) \big\}.
	\end{split}
	\]
	(Note that, by the definition of $\Psf(t)$, the vertex $w$ above is always unique.)
	\item \textbf{Stop.} Terminate the process when $\Asf (t) = \varnothing$. Set $T =t$.
\end{itemize}
\end{defi}

\medskip

Now we are ready to complete the proof of the lower bound on $p_c^{\sss (d)}(\theta )$ in 
Theorem \ref{thm-pc}:
 
\proof[Proof of the lower bound in Theorem \ref{thm-pc}]
Note that $\Dsf(T) = |\cluster(v)|$. Moreover, the line-wise exploration process can be naturally coupled 
to the (modified) Galton-Watson process where the offspring distribution is given by the law of the 
sum of the sizes of $d-1$ independent ERRG clusters minus one (where the $d-1$ accounts for the 
fact that the line we are in has already been explored, and the minus one accounts for the fact that the 
vertex $v^t$ has already been counted). 
This Galton-Watson process has a modification at the root, 
where we consider $d$ independent ERRG cluster sizes, to account for the fact that in the first step 
we have not yet explored any lines. Let $Z_p$ denote the total progeny of this GW-process. 
\ch{A standard argument tells us that $|\cluster(v)|$ is stochastically dominated by $Z_p$, because $Z_p$ ``ignores'' the loops of $|\cluster(v)|$ that do not occur within a line.}
The offspring distribution of this GW-process has mean $\mu=(d-1) (\chil(p)-1)$ (except at the root, where it 
has mean $\mu_{\rho}=d(\chil (p) -1)$), so 
	\eqan{
	\Ep[|\cluster(v)|] &\le \E[Z_p] = 1+\frac{\mu_{\rho}}{1-\mu}=\chil(p) \, \frac{1}{1-(d-1)(\chil(p)-1)}.
	}
Consequently, the value $p_\ell$ that solves the equation
\begin{equation}\label{e:pteqn}
	\theta V^{1/3} = \chil(p) \,  \dfrac{1}{1- (d-1)(\chil(p_\ell)-1)}
\end{equation}
is a lower bound on $p_c^{\sss (d)} (\theta)$.

We insert \eqref{e:chilpsimple} in \eqref{e:pteqn}, to obtain
\begin{equation}
	\theta V^{1/3} = \chil(p) \, \dfrac{1}{1- (d-1)\Big(\dfrac{1}{1-p_\ell(n-1)}
	\Big(1 - \dfrac{2d^2-1}{2(d-1)^3}m^{-1} + O(m^{-2})\Big)-1\Big)}.
\end{equation}
Solving this with respect to $p_\ell$, we find
\begin{equation}
	p_\ell= m^{-1}+\dfrac{2d^2-1}{2(d-1)^2}m^{-2} +O(m^{-1}V^{-1/3}+m^{-3}),
\end{equation}
which gives us the desired lower bound on $p_c^{\sss (d)} (\theta )$. \qed

\medskip
We have found a lower bound on $p_c^{\sss (d)} (\theta )$ by rather explicitly using the product 
structure of the Hamming graph to find a good branching process domination. To find an upper 
bound on $p_c^{\sss (d)} (\theta)$ with the same method would be much more involved, since 
then we would need to thin the GW tree further to take into account the loops outside lines as 
well. Instead, we apply the lace expansion, which is a less direct, but much more robust method.


\section{Bounds on connection probabilities}\label{sect-conn} 

In Section~\ref{sect-up} we show how the lace expansion can be used to express an upper 
bound on $p_c^{\sss (d)}(\theta)$ in terms of products and sums of connection probabilities 
$\P_p(x \leftrightarrow y)$. \ch{In this section, we collect some preliminaries that will be used throughout the analysis in 
Section~\ref{sect-up}. This section is organized as follows.}
In Section~\ref{ssect-tkp} we derive estimates on connection 
probabilities. In Section~\ref{ssect-btpd} we estimate bubble, triangle and polygon diagrams.


\subsection{Connection probabilities}
\label{ssect-tkp}

Given an event $A$ for percolation on $H(d,n)$, we define the event $\{A $ within one line$\}$ to 
be the subset of all configurations $\omega \in A$ such that it can be verified that $\omega \in A$ 
by fixing a line of $H(d,n)$ and then only inspecting the status of the edges of $\omega$ within 
that line. Likewise, we define $\{A$ through multiple lines$\} := A \setminus \{A $ within one line$\}$.
We define $\{v \arf{ \le r} w\}$ to be the event that $w \in \cluster(v)$ and there exists a path of length 
at most $r$. We further write $\{v \arf{ > r} w\}$ (respectively, $\{v \arf{ =r} w\}$) for the event that 
$w \in \cluster(v)$ and there is a (not necessarily shortest) simple path of open edges from $v$ to 
$w$ containing more than (respectively, exactly) $r$ edges.
Not requiring minimality \ch{deviates} a 
bit from the common use of this notation, but for all the purposes of the present paper it makes 
little difference and often simplifies proofs. \ch{This is because the event $\{v \arf{ =r} w\}$ is increasing, whereas the event
that the graph distance equals $r$ is not.}

We start with a proposition about the probability of two points being connected by a path that 
is longer than the mixing time of $H(d,n)$. Given a graph $G$, we define the $t$-step 
\emph{non-backtracking random walk} (NBW) on $G$ started at $x$ as the uniform measure on paths 
$(X_1,\dots,X_t)$ such that $X_1 =x$ and $X_i \neq X_{i-2}$ for all $i \in \{3,\dots,t\}$ (i.e., the 
path never backtracks). For two vertices $x,y$ of $G$, we write $\mathbf{p}_{\sss \mathrm{NBW}}^t(x,y)$ 
for the probability that a $t$-step non-backtracking random walk started at $x$ ends at $y$. Given 
a connected aperiodic graph $G$ and $\alpha \in (0,1)$, we define the \emph{uniform non-backtracking 
mixing time} as
\begin{equation}\label{e:mnotdef}
	 \tmix(G;\alpha) := \min \Big \{ t \colon \max_{x,y} \mathbf{p}_{\sss \mathrm{NBW}}^t(x,y)  
	 \leq (1+\alpha)V^{-1} \Big \} \, .
\end{equation}
In the remainder of this paper we will use 
\begin{equation}\label{e:tmixbd}
	\tmix := \tmix(H(d,n); n^{-1}) = O(\log n),
\end{equation}
so $\alpha = n^{-1}$. The above bound is proved by Fitzner and van der Hofstad  \cite{FitHof13}. For this choice of $\alpha$, the 
following proposition is a direct consequence of \cite[Theorem~1.4 and Lemma~3.14]{HofNac12}:

\begin{prop}[A uniform connection bound]\label{long}
	Consider percolation on $H(d,n)$ with $d \ge 2$ and $p \le p_c^{\sss (d)}$. Then
	\begin{equation}\label{e:unif}
		\Pp \big(x \arfl{> \tmix} y \big) \le \frac{ \chi(p)}{V} (1+ O(m^{-1})).
	\end{equation}
\end{prop}

\noindent
Heuristically this proposition can be explained with the idea that percolation paths in sufficiently 
high-dimensional graphs at criticality look like random walk paths, so if the path is longer than 
the mixing time, then the connection probabilities become uniform over the graph.
\ch{For Hamming graphs, there is little difference between non-backtracking walk and simple random walk, 
so in many of our bounds we use simple random walk instead.}

We proceed with a useful bound on the two-point function:

\begin{prop}[Asymptotics for the two-point function on the Hamming graph]
\label{twopoint}
Consider percolation on $H(d,n)$ with $d \ge 2$ and $p = m^{-1} (1+O(m^{-1}))$ and $p \le 
p_c^{\sss (d)}$. For any $v,w \in \Vcal$,
\begin{equation}
	\tau_p(v-w) = \delta_{v,w} + \frac{d}{d-1} \frac{1}{m} \indi_{\{v,w\} \in \Ecal} 
	+ \frac{\chi(p)}{V}(1+O(m^{-1})) + O(m^{-(d(v,w)\vee 2)}),
\end{equation}
where $d(v,w)$ is the graph distance between $v$ and $w$ on $H(d,n)$.
\end{prop}

\proof
If $v=w$ then $\tau_p(v-w) = \P_p(v \leftrightarrow w) = 1$ by definition. This gives rise to the 
factor $\delta_{v,w}$ above. Let us henceforth assume that $v \neq w$. We divide the event 
$\{v \leftrightarrow w\}$ into three disjoint events as follows:
\begin{itemize}
	\item $A= \big\{v \arfl{> \tmix} w \big\}$,
	\item $B= \big\{v \arfl{\sss\leq \tmix} w \text{ through different lines}\big\}$,
	\item $C= \big\{v \arfl{\sss\leq \tmix} w \text{ within one line} \big\}$.
\end{itemize}
We bound their contributions separately.

By Proposition \ref{long} we have
 \begin{equation}\label{e:PpAbd}
	\P_p(A) = \dfrac{\chi(p)}{V}(1+O(m^{-1})).
\end{equation}
For the bound on $\Pp(B)$ we distinguish three different cases: $d(v,w) =1$, $1<d(v,w) < d$, 
$d(v,w)= d$.

\medskip
\paragraph{\bf Case $1<d(v,w)< d$:}
Write $\mathbf{p}^t(v,w)$ for the probability that a simple random walk started at $v$ is at $w$ after 
$t$ steps, and let $\Pcal_k(v,w)$ denote the set of all simple paths of length $k$ from $v$ to $w$ in 
$H(d,n)$. There are at most $m^k \mathbf{p}^k(u,v)$ such paths, and so 
\begin{equation}\label{e:pathcounter}
	\begin{split}
		\Pp(B) & \le \sum_{k=d(v,w)}^{\tmix} \sum_{\eta \in \Pcal_k(v,w)} 
		\Pp(\eta \text{ is open}) \le  \sum_{k=d(v,w)}^{\tmix} p^k m^k \mathbf{p}^k(v,w)\\
		&= (1+o(1)) \sum_{k=d(v,w)}^{\tmix} \mathbf{p}^k(v,w),
	\end{split}
\end{equation}
where for the last bound we use that $(mp)^k = 1 +o(1)$ for all $k \le \tmix$ by \eqref{e:critwindowBCHHS} 
and \eqref{e:tmixbd}. Define the set 
\begin{equation}\label{e:hyperplane}
	\Hcal(v,w):= \big\{u \in \Vcal : u_i=v_i \text{ for all } i \in \{1,\dots,d\}\text{ such that } v_i=w_i \big\}.
\end{equation}
We can view $\Hcal(v,w)$ as the ``lowest-dimensional hyperplane'' that contains both $v$ and $w$. We 
write
\begin{equation}
	\mathbf{p}^{k}(v,w)=:\mathbf{p}_{\Hcal(v,w)}^{k}(v,w ) 
	+ \mathbf{p}_{\neg \Hcal(v,w)}^{k}(v,w ),
\end{equation}
where $\mathbf{p}_{\Hcal(v,w)}^{k}(v,w )$ is the probability to go from $v$ to $w$ in $k$ steps without 
leaving $\Hcal(v,w)$.

If the walker started from $v$ is to reach $w$, then it will need to take a step in each direction such 
that $v_i \neq w_i$. At each step of a walk on the Hamming graph the probability that the walk stays 
in $\Hcal(v,w)$ and gets closer to $w = (w_1,\dots,w_d)$ is at most $\frac{d(v,w)}{m}$, since it has 
to move in one of at most $d(v,w)$ directions, say direction $j$, exactly to the unique neighbour that 
has $j$th coordinate $w_j$. There are at most $k^{d(v,w)}$ orders in which the distance-decreasing \ch{steps can occur among $k$ steps.}
Therefore
\begin{equation}\label{inL}
	\begin{split}
		\sum_{k=d(v,w)}^{\tmix} \mathbf{p}_{\Hcal(v,w)}^{k}(v,w )
		& \le \sum_{k= d(v,w)}^{\tmix} k^{d(v,w)} \left(\frac{d(v,w)}{m}\right)^{d(v,w)}
		\left(\dfrac{d(v,w)}{d}\right)^{k-d(v,w)}\\
		&\leq \left(\frac{d}{m}\right)^{d(v,w)} \sum_{k= d(v,w)}^{\tmix} k^{d(v,w)} 
		\left(\dfrac{d(v,w)}{d}\right)^{k} = O(m^{-d(v,w)}).
	\end{split}
\end{equation}
If the path from $v$ to $w$ leaves $\Hcal(v,w)$, then it will have to take at least $d(v,w)+1$ steps in 
the direction of $w$. Since $k \leq \tmix$, there are at most $\tmix$ places along the walk where these 
steps can occur, so we can bound
\begin{equation}\label{e:pathcountoutofL}
	\begin{split}
	\sum_{k=d(v,w)}^{\tmix} \mathbf{p}_{\neg \Hcal(v,w)}^{k}(v,w )
	&\leq \tmix \cdot \tmix^{d(v,w)+1} \left( \frac{d }{m}\right)^{d(v,w)+1} \\
	&= O \big( m^{-d(v,w)-1} \tmix^{d(v,w)+2} \big).
	\end{split}
\end{equation}

\medskip
\paragraph{\bf Case $d(v,w) =d$:}
The walk needs to take at least $d$ steps in the direction of $w$. Hence, using the 
same argument as for \eqref{e:pathcountoutofL}, we obtain
\begin{equation}\label{e:pathdimd}
	\sum_{k=d}^{\tmix} (mp)^k \mathbf{p}^{k}(v,w)=O(m^{-d} \tmix^{d+1}).
\end{equation}

\medskip
\paragraph{\bf Case $d(v,w) =1$:}
Given a random walk starting from $v$, we write $Z$ for the (random) first vertex on the walk such that 
$d (Z,w)=2$ and $T$ for the number of steps it took to go from $v$ to $Z$. A standard path-counting 
argument gives an upper bound on $\Pp(B)$ for the case $d(v,w)=1$ as follows:
\begin{equation}\label{out}
	\P_p(B) \leq \sum_{t=1}^{\tmix}\sum_{z\in \Vcal} 
	\P_p (Z= z,  T=t)\sum_{k=1}^{\tmix -t} (mp)^{k+t} \mathbf{p}^{k}(z,w).
\end{equation}
By definition, the vertex $Z$ is unique and $\P_p (Z=z)\neq 0$ only if $d(z,w)=2$, and, by the transitivity 
of $H(d,n)$, for every $k$ we know that $\mathbf{p}^{k}(z,w) = \mathbf{p}^{k}(y,w)$ if $d(z,w)=d(y,w)$. 
Thus
\begin{equation}\label{e:PpBbd}
	\P_p(B) \leq \sum_{t=1}^{\tmix} \P_p (T=t) 
	\sum_{k=1}^{\tmix} (mp)^{\tmix}  \mathbf{p}^{k}(z,w) 
	= O( m^{-(2 \vee d(v,w))}+m^{-d}\tmix^{d+1}),
\end{equation}
where we \ch{have used the assumption that $p \le m^{-1}(1+O(m^{-1}))$} and \eqref{e:tmixbd} for the second bound. \ch{We further assumed $p \ge m^{-1}$ for the first bound. Note that this can be done without loss of generality.}
\medskip

It remains to bound $\Pp(C)$. We consider the cases $d(v,w)=1$ and $d(v,w) \ge 2$. By definition, 
$\Pp(C) =0$ if $d(v,w) \ge 2$, which takes care of the latter case. If $d(v,w)=1$, then the probability 
to connect within that line is simply the two-point function of an ERRG with $n$ vertices and edge 
probability $p= \frac{1+O(m^{-1})}{d(n-1)}$, which is 
\begin{equation}
	\frac{\chi_G(p)-1}{n-1}=\frac{1+O(m^{-1})}{(d-1)(n-1)},
\end{equation}
for every pair of distinct vertices, due to the symmetry of $K_n$. We therefore obtain
\begin{equation}\label{e:PpCbd}
	\P_p(C)= \dfrac{d}{d-1}\dfrac{1+O(m^{-1})}{m} \indi_{\{v,w\} \in \Ecal}.
\end{equation}
Adding $\delta_{v,w}$ and the three bounds in \eqref{e:PpAbd}, \eqref{e:PpBbd} and \eqref{e:PpCbd} 
for $\Pp(A)$, $\P_p(B)$ and $\P_p(C)$, respectively, completes the proof.
\qed


\subsection{Bubble, triangle and polygon diagrams}
\label{ssect-btpd}

The final estimates of this section involve the so-called bubble, triangle, and polygon diagrams. As 
was already alluded to in the introduction, these diagrams, the triangle diagram in particular, are 
very important quantities in the study of high-dimensional percolation. We start with their definition.

Given an integer $i \ge 2$ and vertices $v ,x_1,\dots, x_{i-1}, w \in \Vcal$, we define the 
\emph{$i$-gon diagrams}
\begin{align}
	C_i^{\sss (0)}(v,x_1,\dots,x_{i-1}, w) &:= \P_p(v \leftrightarrow x_1)\cdots 
	\P_p(x_{i-1} \leftrightarrow w),\\
	C_i^{\sss (1)}(u, x_1, \dots, x_{i-1}, w) &:= \sum_{v \colon \{u,v\} \in \Ecal} 
	p C_i^{\sss (0)}(v, x_1, \dots, x_{i-1}, w),\\
	C^{\sss \leq k}_i(v,x_1 , \dots ,x_{i-1}, w)&:= \sum_{k_1+\cdots+k_{i}\leq k}
	\P_p(v \arfl{=k_1} x_1)\times\cdots\times\P_p(x_{i-1} 
	\arfl{=k_{i}} w),\\
	C^{\sss > k}_i(v ,x_1, \dots ,x_{i-1},w)&:= \sum_{k_1+\cdots+k_{i}>k}
	\P_p(v \arfl{=k_1} x_1)\times \cdots \times\P_p(x_{i-1} 
	\arfl{=k_{i}} w),
\end{align}
where in the case $i=2$ we mean $C_2^{\sss (0)}(v,x_1,y)$, etc. Recall \eqref{e:triangledef} and
observe that
\begin{equation}
	\sum_{x,y \in \Vcal} C_3^{\sss (0)}(v,x,y,w) = \nabla_p(v,w) \le \delta_{v,w}
	+ 10\,\dfrac{\chi(p)^3}{V}+O(m^{-1}),
\end{equation}
and recall that by the definition of $p_c(\theta)$ in \eqref{e:pcdef}, we have $\chi(p_c(\theta))^3
=\theta^3V$. Recall, moreover, the following useful bound by Borgs et al.\ in  \cite[Proposition~1.2 and (5.106)]{BorChaHofSlaSpe05b} (adapted to our setting): For all $v, w \in \Vcal$, $p = m^{-1} (1+O(m^{-1} 
+ V^{-1/3}))$, and all $p \le p_c^{\sss (d)}(\theta)$ with $\theta$ sufficiently small such that $\theta^3 
\le \beta_0$ for $\beta_0$ such that \eqref{e:triangledef} holds,
\begin{equation}\label{e:opentri}
	\sum_{x_1, x_2} C_3^{\sss (1)}(v,x_1, x_2, w) = 3\, \frac{\chi(p)^3}{V} + O(m^{-1}).
\end{equation}

We proceed with a bound on ``short'' polygons:

\begin{lem}\label{short}
Consider percolation on $H(d,n)$ with $p= m^{-1} (1+O(m^{-1}))$. Then for all $v,w \in \Vcal$ and 
each integer $i \ge 2$,
\begin{equation}
	\sum_{x_1,\dots,x_{i-1} \in \Vcal}  C_i^{\sss\leq \tmix} (v,x_1,\dots,x_{i-1},w) 
	= \delta_{v,w} +  O (m^{-d(v,w) \vee 1}+ m^{-d} \tmix^{d+i}).
\end{equation}
\end{lem}

\proof
We use the same path counting argument as in the proof of Proposition~\ref{twopoint}. Instead 
of fixing $x_1, \dots, x_{i-1}$ and summing over all possible probabilities, we fix a random walk 
path of length $k \le \tmix$ between $v$ and $w$ and count the possible ways in which that path 
could appear in the above sum. There are at most $(k+1)^{i-1}$ possible ways to mark the path 
with the vertices $x_1,\dots,x_{i-1}$. We thus bound
\begin{equation}
	\sum_{x_1,\dots,x_{i-1} \in \Vcal}  C_i^{\sss\leq \tmix} (v,x_1,\dots,x_{i-1},w) 
	\leq \sum_{k=0}^{\tmix} (k+1)^{i-1}(mp)^k  \mathbf{p}^{k}(v,w). 
\end{equation}
Compare this bound with \eqref{e:pathcounter}. The only difference is \ch{the factor $(k+1)^{i-1}$.} 

\medskip
\paragraph{{\bf Case $d(v,w) < d$:}}
Recall the definition of $\Hcal(v,w)$ in \eqref{e:hyperplane}. We again consider the contributions 
from walks that stay within $\Hcal(v,w)$ and those that do not separately, starting with the \ch{contribution from the walks that}
remain within $\Hcal(v,w)$.

If $v=w$ then $\Hcal(v,w) = \{v\}$, so the only contribution comes from the trivial path $v=x_1=\dots
= x_{i-1}=w$, which gives the term $\delta_{v,w}$. If $1< d(v,w) <d$, then to go from $v$ to $w$ 
the walk needs to take at least one step in each direction $j$ such that $v_j \neq w_j$ and move 
exactly to a neighbour with $j$th coordinate $w_j$. Using a similar argument as in \eqref{inL}, we 
obtain
\begin{equation}
	\sum_{k=d(v,w)}^{\tmix} (k+1)^{i-1}(mp)^k\mathbf{p}_{\Hcal(v,w)}^{k}(v,w )
	=O(m^{-d(v,w)}).
\end{equation}
Note, in particular, that the
\ch{extra factor $(k+1)^{i-1}$}
compared to \eqref{inL} \ch{does} not affect the 
convergence of the sum over $k$:
\ch{it only changes the constant} in the $O(m^{-d(v,w)})$ term.

Next, consider the contribution of walks that leave $\Hcal(v,w)$. If the walk leaves $\Hcal(v,w)$ along 
the path from $v$ to $w$, then it needs to take at least $d(v,w)+1$ distance-decreasing steps. By 
the same argument as \eqref{e:pathcountoutofL}, we thus bound
\begin{equation}
\begin{split}
	\sum_{k=d(v,w)}^{\tmix}(k+1)^{i-1} (mp)^k \mathbf{p}_{\neg \Hcal(v,w)}^{k}(v,w )
	&\leq (1+o(1))\tmix^{d(v,w)+2+i-1}  \left(\frac{d}{m}\right)^{-d(v,w)-1}\\
	&= O(\tmix^{d(v,w)+i+1}  m^{-d(v,w)-1}) 
	= O(m^{-d(v,w)}).
\end{split}
\end{equation}

\medskip
\paragraph{{\bf Case $d(v,w)=d$:}} The walk needs to take at least $d$ steps towards 
$w$ to reach it. By the same argument as for the previous bound, we have
\begin{equation}
	\sum_{k=d(v,w)}^{\tmix} (k+1)^{i-1}(mp)^k \mathbf{p}^{k}(v,w)=O(m^{-d} \tmix^{d+i}).
\end{equation}

Summing these bounds we get the claim. \qed
\medskip

The next lemma combines the above estimates in a convenient form (and is especially 
useful when $i=2$).

\begin{lem}\label{lem:poly}
Consider percolation on $H(d,n)$ with $p \le p_c^{\sss (d)}(\theta)$. For all $i \ge 2$, $j \in \{0,1\}$, 
and $v,w \in \Vcal$,
\begin{equation}
	\begin{split}
		\sum_{x_1, \dots, x_{i-1} \in \Vcal} C^{\sss (j)}_i (v,x_1, \dots, x_{i-1},w) 
		& = \delta_{v,w}\delta_{j,0}  + \theta^{i} V^{i/3 -1} (1 + O(m^{-1}))\\
		&\qquad \quad + O(m^{-d(v,w) \vee 1} + m^{-d} \tmix^{d+i}).
	\end{split}
\end{equation}
\end{lem}

\proof For $j=0$, the lemma follows after combining Proposition~\ref{long} and Lemma~\ref{short} 
with \eqref{e:pcdef}. For $j=1$ we further observe that 
\begin{equation}
	\sum_{u: \{v,u\} \in \Ecal} p = m p = 1 + O(m^{-1}).
\end{equation}
Moreover, the factor $\delta_{v,w}$ does not arise, because this is due to the trivial path $v= x_1 
= \dots=x_{i-1} = w$, which by Proposition~\ref{twopoint} now contributes
\begin{equation}
	C_i^{\sss (1)}(v,v,\dots,v,v) = p\sum_{u : \{v,u\} \in \Ecal} C_3^{\sss (0)}(u,v,\dots,v,v) 
	= p \sum_{u : \{u,v\} \in \Ecal} \tau_p(v-u) = O(m^{-1}).
\end{equation}
\qed

\medskip
Lastly, we derive an improved bound on the triangle diagram in the case where the two intermediate 
points of the triangle are neighbours in $H(d,n)$ and in the case where one intermediate point of the triangle is 
constrained to be a neighbour of a fixed auxiliary point:

\begin{lem}\label{lem:distone}
Consider percolation on $H(d,n)$ with $p \le p_c^{\sss (d)}(\theta)$. For all $v,w \in \Vcal$ and 
$j \in \{0,1\}$,
\begin{equation}\label{e:C3}
	\sum_{x,y : \{x,y\} \in \Ecal} C_3^{\sss (j)} (v,x,y,w)  = O(m^{-1} + V^{-1/3}) 
\end{equation}
and
\begin{equation}\label{e:C3z}
	\sup_z \sum_{x, y : \{y,z\} \in \Ecal} C_3^{\sss (j)}(v,x,y,w) = O(m^{-1} + V^{-1/3}).
\end{equation}
\end{lem}

\proof
We prove \eqref{e:C3} for the case $j=0$. The proof of \eqref{e:C3} for the case where $j=1$ 
and of \eqref{e:C3z} are almost identical, so we leave them to the reader.

Consider first the contribution due the cases where all the connections are due to short 
paths (i.e., shorter than $\tmix$ in total). By Lemma~\ref{short}, this is bounded by $O(m^{-1})$, 
even without using the constraint that $\{x,y\} \in \Ecal$. Next, consider the contribution to the 
left-hand side of \eqref{e:C3} from the case where the path from $y$ to $w$ is longer than $\tmix$,
 i.e.,
\begin{equation}
	\sum_x \P_p(v \leftrightarrow x) \sum_{y : \{x,y\} \in \Ecal} \P_p(x \leftrightarrow y)\, 
	\P_p(y \arfl{> \tmix} w).
\end{equation}
Applying Proposition~\ref{long} to the last term, Proposition~\ref{twopoint} to the middle term, and 
summing over $x$ and $y$, we obtain the upper bound
\begin{equation}\label{e:distone1stbd}
	\chi(p) \, m O(m^{-1}) \, O \Big(\frac{\chi(p)}{V}\Big) = O(V^{-1/3}),
\end{equation}
where the final bound is due to \eqref{e:pcdef}. The contribution due to the case where the path from
$v$ to $x$ is longer than $\tmix$ is the same by symmetry.

To bound the contributions due to the case where the path from $x$ to $y$ is longer than $\tmix$, 
we consider the cases $y=w$ and $y \neq w$ separately. The contribution due to a long path between 
$x$ and $y$ and $y=w$ is given by
\begin{equation}
	\sum_x \P_p(v \leftrightarrow x)\, \P_p(x \arfl{> \tmix} w) 
	= \chi(p) \, O\Big(\frac{\chi(p)}{V}\Big) = O(V^{-1/3}),
\end{equation}
where the bound follows from Proposition~\ref{long} and \eqref{e:pcdef}. The contribution due to a 
long path between $x$ and $y$ and $y \neq w$ is given by
\begin{equation}\label{e:distonelastbd}
	\sum_x \P_p(v \leftrightarrow x) \sum_{\substack{y: \{x,y\} \in \Ecal,\\ y \neq w}} 
	\P_p(x \arfl{> \tmix} y)\, \P_p(y \leftrightarrow w) 
	= \chi(p) \, m O\Big(\frac{\chi(p)}{V}\Big)\,O(m^{-1}) = O(V^{-1/3}),
\end{equation}
where the bound again follows from Propositions~\ref{long} and \ref{twopoint}, and \eqref{e:pcdef}. 

Adding the bound $O(m^{-1})$ due to short paths and the bounds 
\eqref{e:distone1stbd}--\eqref{e:distonelastbd}, we complete the proof of \eqref{e:C3} for the case 
$j=0$.\qed


\section{The upper bound on $p_c (\theta)$ via the lace expansion}\label{sect-up}

In Section~\ref{ssect-bg} we recall the background of the lace-expansion technique and
state a proposition in which we estimate lace-expansion coefficients. In Section~\ref{ssect-ub} 
we use this proposition to prove the upper bound on $p_c (\theta)$. The proof of the proposition 
is given in Sections~\ref{ssect-est1}--\ref{subsect:Pin}.


\subsection{Background}
\label{ssect-bg}

The lace expansion is a method originated by Brydges and Spencer to study self-avoiding walk 
\cite{BrySpe85}, and was first applied to percolation by Hara and Slade \cite{HarSla90}. Hara and 
Slade's method gives an expansion for the percolation two-point function $\tau_p(x)$. The version 
of the lace expansion that we use here was derived by Borgs et al.\ \cite{BorChaHofSlaSpe05b}, where it is proved 
for the Hamming graph, among others, that for any $p \in [0,1]$,
\begin{equation}\label{e:LE}
	\tau_p (x) = \delta_{0,x} + m (J_p * \tau_p)(x) + m (\Pi_p * J_p * \tau_p)(x) + \Pi_p(x),
\end{equation}
where $(f * g)(x) = \sum_{y \in \Vcal} f(y) g(x-y)$ denotes the convolution between $f$ and $g$, $J_p(x-y) := \Pp(\{x,y\}$ is open$)$, 
and $\Pi_p(x)$ 
is the so-called \emph{irreducible two-point function} or \emph{lace-expansion coefficient}. The lace expansion further determines that
$\Pi_p(x)$ is given by the alternating series
\begin{equation}\label{e:Piseries}
	\Pi_p(x) = \sum_{N=0}^\infty (-1)^N \Pi_p^{\sss (N)}(x),
\end{equation}
and the \emph{lace-expansion coefficients} $\Pi_p^{\sss (N)}$ have a well-defined structure that will 
play an important role in the determination of the upper bound we derive here. Define the \emph{discrete 
Fourier transform} of a function $f : \Vcal \to \mathbb{R}$ as $\hat f (k) = \sum_{x \in \Vcal} \e^{i k \cdot x} 
f(x)$ for $k \in \mathbb{R}^d$. Taking the Fourier transform of \eqref{e:LE}, we obtain
\begin{equation}
	\hat \tau_p(k) = \frac{1 + \hat \Pi_p(k)}{1 - m \hat J_p(k) (1 + \hat \Pi_p(k))}.
\end{equation}
Observe that $\hat \tau_p(0) = \sum_{x} \tau_p(x) = \chi(p)$ and $\hat J_p(0) = mp$, so setting $k \equiv 0$ 
and $p \equiv p_c^{\sss (d)}(\theta)$, and applying \eqref{e:pcdef}, we obtain
\begin{equation}\label{e:pcexp}
	m p_c(\theta) = \frac{1}{1+ \hat \Pi_{p_c}(0)}  + \theta^{-1} V^{-1/3}.
\end{equation}
(We will henceforth only consider $\hat \Pi_p (k)$ at $k=0$ and therefore will not write the argument 
anymore.)
Combining \eqref{e:Piseries} and \eqref{e:pcexp} with the following proposition allows us to determine 
the upper bound on $p_c^{\sss (d)}(\theta)$ for $\theta$ sufficiently small:

\begin{prop}[Bounds on the lace-expansion coefficients]\label{prop:Pi}
Consider percolation on $H(d,n)$ with $d \ge 3$. Let $\theta$ be such that $3 \theta^3 (1+ 10 \, \theta^3) 
< 1$ and $\theta^3 \le \beta_0$ for $\beta_0$ such that \eqref{e:triangledef} holds, and let $p$ be such 
that $p = m^{-1} (1+O(m^{-1}+V^{-1/3}))$ and $p \le p_c^{\sss (d)}(\theta)$. Then 
\begin{align}
	\hat\Pi^{\sss (0)}_p  &\geq\, \dfrac{2d-1}{2(d-1)^2}m^{-1} +O(m^{-2} + V^{-1/3}), \label{pi0}\\
	\hat\Pi^{\sss (1)}_p  &\leq\,  \dfrac{d^2+d-1}{(d-1)^2}m^{-1} +O(m^{-2} + V^{-1/3}),\label{pi1} \\
	\sum_{N=2}^\infty \hat\Pi^{\sss (N)}_p  &=\, O(m^{-2} + V^{-1/3}). \label{pi2}
\end{align}
\end{prop}
Note that we do not consider $d=2$ in the above proposition. This is of no consequence, as our main 
theorem is already proved for $d=2,3$ in \cite{BorChaHofSlaSpe05a}. The proofs below do apply to 
the case $d=2$, but the bounds become slightly less sharp (because Lemma~\ref{lem:poly} is less 
sharp when $d=2$.)


\subsection{Proof of the upper bound on the critical point}
\label{ssect-ub}

Before starting with the proof of Proposition \ref{prop:Pi}, which constitutes the bulk of what remains 
of this paper, let us complete the proof of the upper bound in Theorem \ref{thm-pc}:

\proof[Proof of the upper bound in Theorem \ref{thm-pc} subject to Proposition \ref{prop:Pi}]
As was mentioned above, Theorem~\ref{thm-pc} is already proved for $d =2,3$ in 
\cite{BorChaHofSlaSpe05a}, so let $d \ge 4$. Suppose first that $\theta$ is such that 
$3 \theta^3 (1+ 10 \, \theta^3) < 1$ and $\theta^3 \leq \beta_0$ for $\beta_0$ such that 
\eqref{e:triangledef} holds. Using \eqref{e:Piseries} and Proposition~\ref{prop:Pi}, we bound 
\begin{equation}
	\begin{split}
		\hat\Pi_{p_c^{\sss (d)}} &= \hat\Pi^{\sss (0)}_{p_c^{\sss (d)}} 
		- \hat\Pi^{\sss (1)}_{p_c^{\sss (d)}}  + O(m^{-2} + V^{-1/3}) \\
		& \geq \dfrac{2d-1}{2(d-1)^2}m^{-1} - 
		\dfrac{d^2+d-1}{(d-1)^2}m^{-1}+O(m^{-2} + V^{-1/3})\\
		& \geq - \dfrac{2d^2-1}{2(d-1)^2}m^{-1}+O(m^{-2} + V^{-1/3}).
	\end{split}
\end{equation}
Applying this bound to \eqref{e:pcexp}, we get
\begin{equation}
	\begin{split}
		p_c^{\sss (d)} (\theta) &\le \frac{m^{-1}}{1-\frac{2d^2-1}{2(d-1)^2}m^{-1} 
		+O(m^{-2} + V^{-1/3})}+ \theta^{-1} m^{-1} V^{-1/3}\\
		& = m^{-1}+\frac{2d^2-1}{2(d-1)^2}m^{-2}+O( m^{-1} V^{-1/3}+m^{-3}),
	\end{split}
\end{equation}
which gives the upper bound for $\theta$ such that $3 \theta^3 (1+ 10 \, \theta^3) < 1$ and 
$\theta^3 \leq \beta_0$. 

We have determined $p_c (\theta )$ for some $\theta \in (0, \infty )$. Borgs et al.\  \cite[Theorem 1.1]{BorChaHofSlaSpe05a} show that $p_c (\theta ')=p_c(\theta ) (1+O(V^{-1/3}))$ \ch{for every $\theta ' \in (0,\infty )$,} and so we obtain the desired upper bound for arbitrary $\theta$, which completes the proof of the upper bound. \qed

\medskip
The remainder of this section is devoted to the proof of Proposition~\ref{prop:Pi}, \ch{which is divided into}
three further subsections.


\subsection{Analysis of $\hat\Pi^{\sss{( 0)}}_p$: proof of \eqref{pi0}.} 
\label{ssect-est1}

The lace-expansion method yields (see \cite[Section~3.2]{BorChaHofSlaSpe05b})
\begin{equation}\label{pi0d}
 	\hat\Pi^{\sss (0)}_p := \sum_{x\neq 0} \P_p (0 \Longleftrightarrow x),
\end{equation}
where $v \Longleftrightarrow w$ denotes the event that there exist two edge-disjoint paths between 
$v$ and $w$. In this case we say that $v$ and $w$ are \emph{doubly connected.} This is equivalent 
to the event that there exists an edge-disjoint cycle of open edges containing both $v$ and $w$. 

Given a graph $G$ and a vertex $x$ in $G$, write $\Lcal_k(G;x)$ for the set of vertex-disjoint cycles 
of length $k$ in $G$ that contain $x$. Note that $|\Lcal_k(K_n; 0)| = \frac{(n-1)!}{2(n-k)!}$, because 
we must choose the $k-1$ other vertices for the cycle from $n-1$ possible choices, and although their 
order matters, the direction of the cycle does not, which explains the factor $\frac12$. By only considering
vertex-disjoint cycles in $\hat \Pi_p^{\sss (0)}$ that are contained within a line that intersects $0$, we 
can apply a similar argument to that of \eqref{e:Splb} to bound
\begin{align}
	\hat \Pi_p^{\sss (0)} &\ge \sum_{x \neq 0} \P_p (0 \Longleftrightarrow x 
	\text{ within one line})\\
	&\ge d \sum_{k=3}^n (k-1) \sum_{L_k \in \Lcal_k(K_n ; 0)} 
	\bigg[\P_p(L_k \in \Lcal_k(\cluster(0) ; 0))\nn\\
	& \qquad\qquad\qquad\qquad\qquad\qquad - \P_p\Big(\bigcup_{\{x,y\} \subset \Vcal(L_k)} 
	\{L_k \in \Lcal_k(\cluster(0) ; 0)\} \circ \{x \leftrightarrow y\}\Big) \bigg]\nn\\
	& \ge d \sum_{k=3}^n (k-1) \frac{(n-1)!}{2(n-k)!} p^k \Big(1 - \tfrac{1}{2} k(k-1) 
	\P_p(1 \in \cluster(0))\Big),\nn
\end{align}
where in the second inequality $\cluster(0)$ is to be viewed as the cluster of vertex $0$ in $G(n,p)$, and 
we use that every cycle of length $k$ that passes through $0$ passes through $(k-1)$ other vertices, 
and in the third inequality $\P_p(1 \in \cluster(0))$ denotes the probability that vertices $0$ and $1$ of 
$K_n$ are in the same cluster in $G(n,p)$. Use that $p = m^{-1}(1+O(m^{-1} + V^{-1/3}))$ by assumption, 
that $\frac{(n-1)!}{(n-k)!} = (n-1)^{k-1} - O(k n^{k-2})$, and that $\P_p(1 \in \cluster(0)) = \frac{d}{(d-1)n} 
(1+ O(m^{-1} + V^{-1/3}))$ by Theorem~\ref{secondorder}, to bound (after a short computation)
\begin{equation}
	\hat \Pi_p^{\sss (0)} = \frac{2d-1}{2(d-1)^2} m^{-1} (1 + O(m^{-1} + V^{-1/3})).
\end{equation}
\qed


\subsection{Analysis of $\hat\Pi^{\sss{(1)}}_p$: proof of \eqref{pi1}.}
\label{ssect:est2}

We start with a few definitions:

\begin{defi}[Connections through a subset and pivotal edges]\color{white}.\color{black}
	\begin{itemize}
		\item 
			We say that two vertices $x,y$ are \emph{connected through} a 
			set $\Wcal \subset \Vcal$, and write $x \overset{\Wcal}{\longleftrightarrow} y$,
			if $x \leftrightarrow y$ occurs and all the paths connecting the two vertices in 
			the percolation configuration have at least one vertex in $\Wcal$.
		\item 
			Given a vertex $x$ and an edge $e$, we define the set $\tilde\cluster^e (x)$ 
			as the percolation cluster of $x$ in the (possibly modified) percolation 
			configuration where the edge $e$ is set to closed. 
		\item 
			Given two vertices $x,y$ and a directed edge $(u,v)$, we say that $(u,v)$ is
			\emph{pivotal} for $x \leftrightarrow y$ if $x \leftrightarrow y$ occurs on the 
			(possibly modified) configuration where $\{u,v\}$ is set to open, 
			and $x \leftrightarrow u$ and $v \leftrightarrow y$ occur, but $x \leftrightarrow y$
			does not occur on the (possibly modified) configuration where $\{u,v\}$ is set to 
			closed. Note that the direction of the edge is important and that
			$\{(u,v)\text{ is pivotal for } x \leftrightarrow y\} \neq \{ ( v,u )\text{ is pivotal for } x 
			\leftrightarrow y \}$.
	\end{itemize}
\end{defi}

We define the event 
\begin{equation}\label{e:Eprime}
	E' (v,x;\Wcal) := \{v \overset{\Wcal}{\longleftrightarrow} x\} \cap \{ \nexists 
	\text{ pivotal } (u',v') \text{ for } v \leftrightarrow x\ \text{ such that } v 
	\overset{\Wcal}{\longleftrightarrow} u' \}.
\end{equation}
This event is a central object of the percolation lace expansion (see e.g.\ \cite[(1.36)]{HarSla90}).
We can now give the definition of $\hat\Pi_p^{\sss (1)}$ from \cite{BorChaHofSlaSpe05b}:
\begin{equation}\label{e:Pinest}
	\hat\Pi_p^{\sss (1)}:= \sum_{x}p\sum_{u,v : \{u,v\} \in \Ecal} \E_{\sss 0} 
	[ \indi_{\{0 \Leftrightarrow u\} } \P_{\sss 1} (E' (v,x;\tilde\cluster_0^{(u,v)}(0))].
\end{equation}
The subscripts $0$ and $1$ indicate that the events happen on two \emph{distinct} percolation 
configurations on the Hamming graph, and that these percolation configurations are ``nested'' in 
such a way that $\tilde \cluster_0^{(u,v)}$ is a set-valued random variable with respect to $\E_0$, 
but for any fixed realization of $\tilde \cluster_0^{(u,v)}$ it is viewed as a deterministic set with respect 
to $\P_1$ (see \cite{BorChaHofSlaSpe05b} for a more in-depth discussion of this construction). 
Our proof strategy will be to first identify and bound the main contributions to $\hat \Pi_p^{\sss (1)}$. 
For this we can apply the BK-inequality to all the events in the formulation of  $\hat\Pi^{\sss{(1)}}_p$ 
and replace the complicated probabilities that appear there by products of two-point functions. 

The analysis of $\hat\Pi_p^{\sss (1)}$ proceeds by reduction of the complicated event to a collection 
of disjointly occurring two-point events, followed by a repeated application of the BK-inequality. Such 
bounds are standard, and are known in the literature as \emph{diagrammatic bounds.} 

\begin{figure}[htbp]
\includegraphics[width=0.40\textwidth]{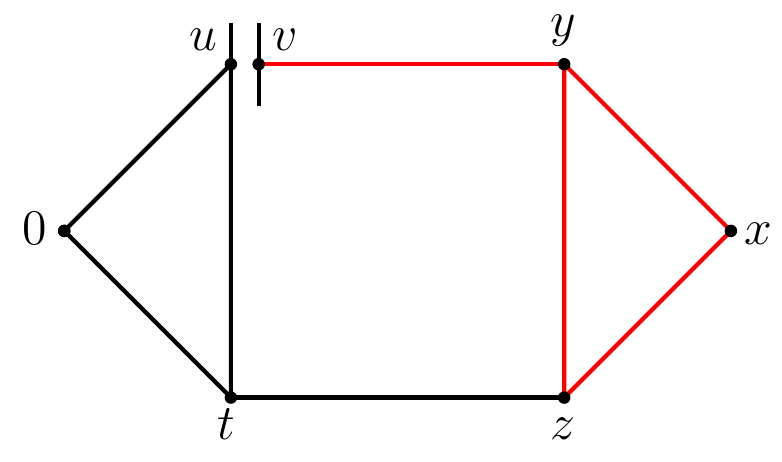}
\caption{Diagrammatic description of the events contributing to $\hat\Pi^{\sss{(1)}}_p$. Black
lines refer to connections that occur in level $0$ of percolation, red lines to connections that 
occur in level $1$. The directed edge $(u,v)$ is represented by two vertical dashes.
\label{diaP1}}
\end{figure}

Inspecting \eqref{e:Eprime} and \eqref{e:Pinest}, we see that on the percolation measure $\P_1$ we 
need that $v$ is connected to $x$ through $\tilde \cluster_0^{(u,v)}(0)$, that any vertex $z \in \tilde 
\cluster^{(u,v)}_0 (0)$ on the path from $v$ to $x$ must be doubly connected to $x$ (otherwise there 
would exist a pivotal $(u', v')$), and that on a path from $z$ to $x$ there must be a vertex $y$ that is 
connected to $v$. Moreover, on the percolation measure $\P_0$, the vertex $z$ is connected to $0$, 
and $0$ is doubly connected to $u$, so there must exist exist a vertex $t$ on a path from $0$ to $u$ 
such that $t$ is connected to $z$. See Figure~\ref{diaP1}. 

The main contribution will come from the \st{simple} case when $0 = u =t$ and $x=y=z$. To get a sharp 
bound on this term, it is necessary to use that \eqref{e:Eprime} implies that the connection from $0$ 
to $x$ on the $\P_0$-percolation configuration is not allowed to use the edge $\{0,v\}$. We will therefore 
bound this contribution by
\begin{equation}
	M:=  p \sum_{x} \sum_{v : \{0,v\} \in \Ecal} \P_p (0 \leftrightarrow x 
	\text{ without using } \{0,v\}) \P_p( v \leftrightarrow x).
\end{equation}
Following the same derivation as in \cite[Section 4.1]{BorChaHofSlaSpe05b}, but isolating the main contribution, we can use the BK-inequality to derive the following upper bound from \eqref{e:Pinest}:
\begin{equation}\label{e:Pionediag}	
	\hat\Pi_p^{\sss (1)} \leq M+  \sum_{u,t,z,y,x}  C_3^{\sss (0)} (0,u,t,0) 
	C_3^{\sss (1)}(u,y,z,t)C_2^{\sss (0)}(y,x,z)
	(1-\delta_{0,u} \delta_{0,t} \delta_{x,y} \delta_{y,z}). 
\end{equation}

Our bound on \eqref{e:Pionediag} will rely heavily on the kind of path-counting methods that we 
applied in the previous sections. 

We proceed by bounding the main term $M$. We claim that
\begin{equation}\label{e:main1}
	\begin{split}
		M &\le p \sum_{v \colon \{0,v\} \in \Ecal} \Big(\big( p + 
		\P_p(v \overset{\ge 2}{\longleftrightarrow} 0)\big) + 
		\P_p(0 \overset{\ge 2}{\longleftrightarrow} x) + 
		\sum_{\substack{x : x \neq v, \\ x\neq 0}}\, \P(0 \leftrightarrow x) 
		\P_p(v \leftrightarrow x) \Big)\\
		& \le p \sum_{v \colon \{0,v\} \in \Ecal}  \Big(p + \sum_x C_2^{\sss\geq 2} (0,x,v)\Big).
	\end{split}
\end{equation}
The first term in the first inequality is due to the case $x=0$, the second term to $x=v$, and the third 
term to the remaining cases. Note that if we had not restricted the connection from $0$ to $x$ to occur 
without using $\{0,v\}$, then the second term would have had an additional $p$, so the same argument 
would have given the upper bound $p \sum_{v}  \big(2p + \sum_x C_2^{\sss \geq 2} (0,x,v)\big)$, which, as 
it will turn out, is not sharp enough.

We use the (by now) familiar path-counting estimates to bound the right-hand side of \eqref{e:main1} 
by
	\eqan{\label{e:M123}
	&p \sum_{v \colon \{0,v\} \in \Ecal}   \Big(p + \sum_x C_2^{\sss \geq 2} (0,x,v)\Big) \\
	 &\quad\le p^2 m + p\shift \sum_{v \colon \{0,v\} \in \Ecal} \shift\Big(\sum_{k =2}^\infty  
	 (k+1)\mathbf{p}_{\Hcal(0,v)}^k(0,v) 
	 + \sum_{k =3}^{\tmix}  (k+1)\mathbf{p}_{\neg \Hcal(0,v)}^k(0,v) 
	 + \sum_x C_2^{\sss > \tmix}(0,x,v) \Big)\nn\\
	  &\quad=: p^2 m + M_1 + M_2 + M_3,\nn
	}
where the factors $(k+1)$ in the second and third term on the right-hand side are due to interchanging 
the sum over $k$ with the sum over $x$. The term $M_1$ is the contribution from walks that are 
constrained to remain within one line:
\begin{equation}
	\begin{split}
		M_1 &\le p m \sum_{k=2}^{\infty} (k+1) p^k (n-1)^{k-1}\\
		& = p \sum_{k=2}^\infty (k+1) d^{-k+1} (1 + O(m^{-1} + V^{-1/3}))\\
		& =  \dfrac{3d -2}{(d-1)^2} m^{-1} (1+O(m^{-1} + V^{-1/3})). 
	\end {split}
\end{equation}
Using that $d(v,w) \ge 1$ and using an estimate similar to \eqref{out}, we further bound
\begin{equation}
	M_2 \leq p\sum_{v : \{0,v\} \in \Ecal }\sum_{k=3}^{\tmix} (mp)^k k(k+1) (n-1)^{-1} 
	= O( m^{-2}).
\end{equation}
Next, we use Proposition~\ref{long} and \eqref{e:pcdef} to bound
\begin{equation}
	\begin{split}
		M_3 &=p\sum_{v : \{0,v\} \in \Ecal}\sum_{x} C_2^{\sss > \tmix}(0,x,v)
		= pm \dfrac{\chi(p)^2}{V}(1+O(m^{-1}))\\
		& \leq \theta^2V^{-1/3} (1+O((m^{-1}+V^{-1/3})).
	\end{split}
\end{equation}
Inserting the bounds for $M_1$, $M_2$ and $M_3$ into \eqref{e:M123}, we conclude that
\begin{equation}\label{e:Mbd}
	M \le p^2 m + \dfrac{3d -2}{(d-1)^2} m^{-1}+ O( m^{-2} + V^{-1/3}) 
	= \dfrac{d^2+d-1}{(d-1)^2}m^{-1} +O(m^{-2} + V^{-1/3}),
\end{equation}
and so it remains to show that the other terms in \eqref{e:Pionediag} are error terms.

We split the remaining contributions on the right-hand side of \eqref{e:Pionediag} according to the 
relative locations of $y$ and $z$ as follows: $y=z$, $d(y,z) =1$, $d(y,z) \ge 2$.

\medskip
\paragraph{{\bf Case $d(y,z)\geq 2$:}} Apply Lemma~\ref{lem:poly} twice and \eqref{e:opentri} 
once to obtain
\begin{equation}
	\sum_{t,u,y,z}C_3^{\sss (0)} (0,u,t,0)  C_3^{\sss (1)}(u,y,z,t)
	\sup_{ z': d(y,z')\geq 2 } \sum_{x} C_2^{\sss (0)} (y,x,z') = O(m^{-2} + V^{-1/3}).
\end{equation}

\medskip
\paragraph{{\bf Case $d(y,z)=1$:}} Apply Lemma~\ref{lem:poly} twice and Lemma~\ref{lem:distone} 
once to obtain
	\eqan{
	&\sum_{t,u}C_3^{\sss (0)} (0,u,t,0)\sup_{t'} \Big( \sum_{y,z: \{y,z\} \in \Ecal} 
	C^{\sss (1)}_3(u,y,z,t') \Big) \sup_{z' : d(y,z')=1}\Big( \sum_{x} C_2^{\sss (0)} (y,x,z')\Big)\nn\\
 	&\qquad\qquad= O( m^{-2} + V^{-1/3}).
	}

\medskip
\paragraph{{\bf Case $y=z$:}} Split again, according to whether $d(t,u) \ge 2$, $d(t,u)=1$, $t=u$. To 
bound the case $y=z$ and $d(t,u)\geq 2$, apply Lemma \ref{lem:poly} three times to obtain
\begin{equation}
	\sum_{t,u}C_3^{\sss (0)}(0,u,t,0) \sup_{t': d( t',u )\geq 2} \Big( \sum_{y}
	C_2^{\sss (1)}(u,y,t')\Big) \sum_{x} C_2^{\sss (0)}(y,x,y) 
	= O(m^{-2} + V^{-1/3}).
\end{equation}
Similarly, to bound the case $y=z$ and $d(t,u)=1$, apply Lemmaa~\ref{lem:poly} twice and 
Lemma~\ref{lem:distone} once to obtain
\begin{multline}
	\sum_{t,u : \{t,u\} \in \Ecal}C_3^{\sss (0)}(0,u,t,0) \sup_{t': d( t',u )=1} 
	\Big( \sum_{y}C_2^{\sss (1)}(u,y,t')\Big) \sum_{x} C_2^{\sss (0)}(y,x,y)
	= O(m^{-2} + V^{-1/3}).
\end{multline}
All that remains is the case $y=z$ and $t=u$. Recall that in $M$ we have already accounted for the 
term where $t=u=0$ and $x=y=z$. Applying Proposition~\ref{twopoint} with the constraint $t=u \neq 0$ 
to the first sum in \eqref{e:Pionediag} and Lemma~\ref{lem:poly} to the second and third sum, we 
obtain
\begin{equation}
	\sum_{t \neq 0} C_2^{\sss (0)}(0,t,0)  \sum_y C^{\sss (1)}_2(t,y,t) \sum_{x} 
	C_2^{\sss (0)}(y,x,y) = O(m^{-2} + V^{-1/3}).
\end{equation}
The case $t=u=0$ and $y=z \neq x$ are analogous and give the same bound.

We have thus shown that all the remaining terms in \eqref{e:Pionediag} are of order $O(m^{-2} + V^{-1/3})$, 
which, combined with the bound \eqref{e:Mbd} on the main term $M$, completes the proof of \eqref{pi1}. 
\qed

\begin{figure}[htbp]
\includegraphics[width=0.9\textwidth]{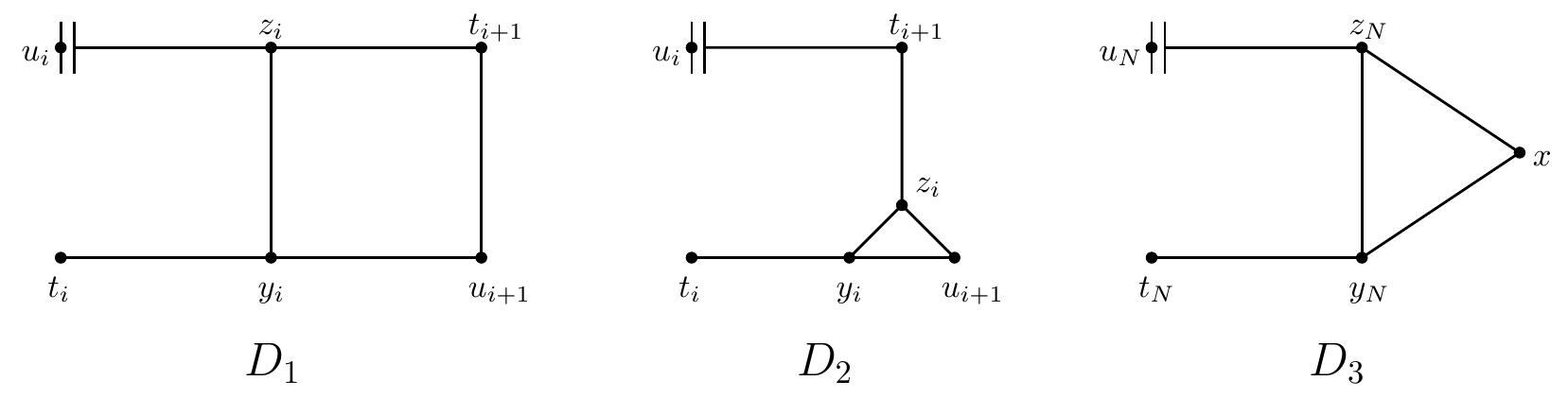}
\caption{Diagrammatic descriptions of $D_1$, $D_2$, and $D_3$. \label{diaPN}}
\end{figure}


\subsection{Analysis of $\hat\Pi^{\sss{(N)}}_p$ for $ N \geq 2$: proof of \eqref{pi2}.}
\label{subsect:Pin}

\begin{figure}[htbp]
\includegraphics[width=0.90\textwidth]{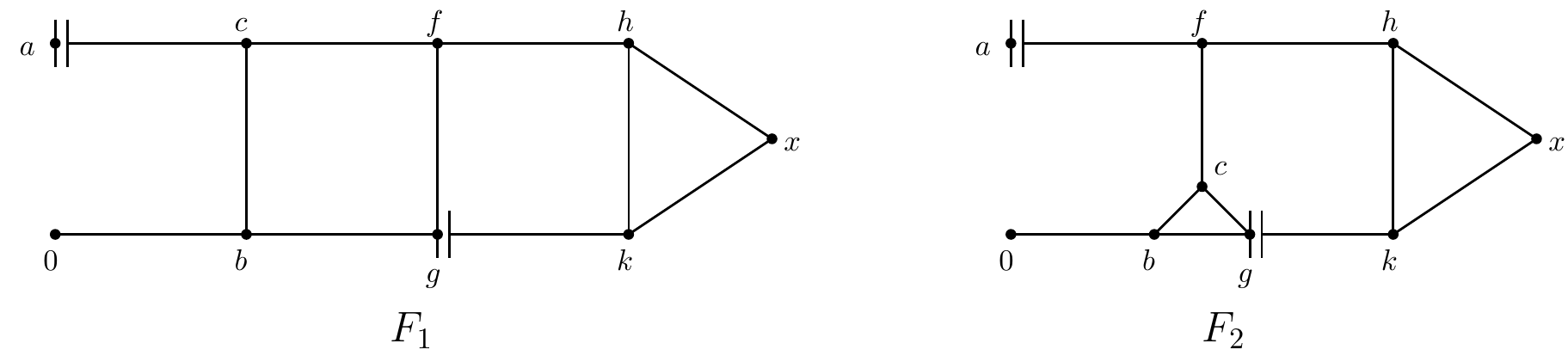}
\caption{Diagrammatic description of the summands in $F_1$, $F_2$. \label{diaPN1}} 
\end{figure}

To investigate $\hat\Pi^{\sss{(N)}}_p$ for $N \geq 2$, we again bound the events in terms of products 
of two-point functions. We define the quantities
\begin{align}
	D_1(t_i,u_i,y_i,z_i,t_{i+1},u_{i+1})
	&:= C_3^{\sss (1)}(t_i,y_i,z_i,u_i) C_3^{\sss (0)} (y_i,t_{i+1},u_{i+1},z_i),\\
	D_2(t_i,u_i,y_i,z_i,t_{i+1},u_{i+1})
	&:=  C_5^{\sss (1)}(t_i,y_i,u_{i+1},z_i,t_{i+1},u_i) \P_p (y_i \leftrightarrow z_i),\\
	D_3(t_N, u_N, y_N,z_N,x)
	&:= C_3^{\sss (1)}(t_N,y_N,z_N,u_N)C_2^{\sss (0)}(y_N,x,z_N).
\end{align}
See Figure~\ref{diaPN}. We write the bounds on $\hat\Pi^{\sss{(N)}}_p$ from 
\cite[Section 4]{BorChaHofSlaSpe05b} in the current notation:
	\eqan{
	\label{e:PiNdiag}
		\hat\Pi_p^{\sss{(N)}}&\leq \sum_{t_1, \dots, t_N} \sum_{u_1, \dots, u_N} 
		\sum_{y_1,\dots, y_N} \sum_{z_1,\dots,z_N} \sum_x  C_3^{\sss (0)} (0,t_1,u_1,0) \\ 
		& \quad \times \prod_{i=2}^{N}[ D_1(t_i,u_i,y_i,z_i,t_{i+1},u_{i+1})
		+D_2(t_i,u_i,y_i,z_i,t_{i+1},u_{i+1}) ] D_3 (t_N,u_N,y_N,z_N,x).\nn
	}
Our strategy for bounding \eqref{e:PiNdiag} will be to first show that the ``tails'' of the diagrams, i.e.,
\begin{align}
	F_1 &:= \sup_{a}\sum_{b,c,f,g,h,k,x}D_1(0,a ,b,c,f,g)D_3(f,g,h,k,x),\\
	F_2 &:= \sup_{a}\sum_{b,c,f,g,h,k,x}D_2(0,a,b,c,f,g)D_3(f,g,h,k,x),
\end{align}
(see Figure~\ref{diaPN1}) are of order $O(m^{-2} + V^{-1/3})$ because of the bounds derived in 
Section~\ref{sect-conn}, and then to bound what remains of \eqref{e:PiNdiag} by $O(\theta^{N-2})$ 
with the help of repeated applications of the bounds in \eqref{e:triangledef} and \eqref{e:opentri}. 
Summing $\hat \Pi^{\sss (N)}_p$ over $N$, we will thus also get $O(m^{-2} + V^{-1/3})$ when 
$\theta$ is sufficiently small.

To bound $F_1$, we treat the cases $h = k$, $d(h,k) = 1$, $d(h,k) \ge 2$ separately.

\medskip
\paragraph{{\bf Case $d(h,k) \ge 2$:}} Apply Lemma~\ref{lem:poly} to the second sum in $D_3$, 
and apply \eqref{e:triangledef} and \eqref{e:opentri} to the remaining three triangle diagrams. This, 
combined with \eqref{e:pcdef}, yields the bound $O(m^{-2} + V^{-1/3})$. 

\medskip
\paragraph{{\bf Case $d(h,k)=1$:}} Apply Lemma~\ref{lem:distone} to the first sum in $D_3$, 
Lemma~\ref{lem:poly} to the second sum in $D_3$, and \eqref{e:triangledef} and \eqref{e:opentri} 
to the sums in $D_1$, to again obtain the bound $O(m^{-2} + V^{-1/3})$.

\medskip
\paragraph{{\bf Case $h=k$:}} Apply Lemma~\ref{lem:poly} to the second sum in $D_3$, to obtain 
the bound
\begin{equation}\label{e:F1gish}
	\sup_a \sum_{b,c,f,g,h} C_3^{\sss (1)}(0,b,c,a) C_3^{\sss (0)}(b,g,h,c) 
	C_2^{\sss (1)}(f,g,h) O(1).
\end{equation}
Compare this with the bound \eqref{e:Pionediag} on $\hat \Pi_p^{\sss (1)}$. The main difference is 
that the placement of the extra open edges is different (i.e., where the zeros and ones are in the 
superscripts). We can analyze \eqref{e:F1gish} essentially in the same way as \eqref{e:Pionediag}. 
The most important difference is that here the term due to $0=a=b=c$ and $f=g=h$ does \emph{not} 
give main contribution, because $j=1$ in the first term here, and so by Lemma~\ref{lem:poly} the term 
in \eqref{e:F1gish} that corresponds to $M$ in \eqref{e:Pionediag} is of order $O(m^{-2} + V^{-1/3})$ 
here. Following the same steps as in the bound on \eqref{e:Pionediag} above, we conclude that all 
other terms are also of order $O(m^{-2} + V^{-1/3})$, and hence,
\begin{equation}\label{e:F1bd}
	F_1 = O(m^{-2} + V^{-1/3}).
\end{equation}

To bound $F_2$ we again treat the cases $h=k$, $d(h,k)=1$, $d(h,k) \ge 2$ separately. 

\medskip
\paragraph{{\bf Cases $d(h,k)=1$ and $d(h,k) \ge 2$:}} The sums in $D_3$ can be bounded in the same way 
as above, to yield a factor $O(m^{-2} + V^{-1/3})$. In \cite[(4.47)]{BorChaHofSlaSpe05b} it is proved 
that
\begin{equation}\label{e:D2bd}
	\sup_a \sum_{b,c,f,g} D_2(0,a,b,c,f,g) \le (1 + 10\, \theta^3 + O(m^{-1}))
	(3 \theta^3 + O(m^{-1})),
\end{equation}
so the bound on $D_2$ is $O(1)$. Consequently, these cases contribute $O(m^{-2} + V^{-1/3})$ to 
$F_2$, as required.

\medskip
\paragraph{{\bf Case $h=k$:}} The final case is more subtle. We consider the contributions from $f=g$, 
$d(f,g)=1$, $d(f,g) \ge 2$ separately. 

If $h=k$ and $d(f,g) \ge 2$, then we may apply Lemma~\ref{lem:poly} 
to the first sum in $D_3$ for a factor $O(m^{-2} + V^{-1/3})$, and to the second sum in $D_3$ for a factor 
$O(1)$. Further applying \eqref{e:D2bd} to the $D_2$ term in $F_2$, we find that this case contributes 
$O(m^{-2} + V^{-1/3})$. 

The contribution due to $h=k$ and $f=g$ is given by
\begin{equation}
	\sup_a \sum_{b,c,f,h,x} C^{\sss (1)}_5(0,b,f,c,f,a) \P_p(b \leftrightarrow c) 
	C_2^{\sss (1)}(f,h,f) C_2^{\sss (0)}(h,x,h).
\end{equation}
By Lemma~\ref{lem:poly}, the third term is bounded by $O(m^{-1} + V^{-1/3})$ and the fourth term is 
bounded by $O(1)$. To bound the two remaining terms, we write
	\eqan{
	&\sup_a \sum_{b,c,f} C^{\sss (1)}_5(0,b,f,c,f,a)\, \P_p(b \leftrightarrow c)\\
	&\qquad= \sup_a \sum_{b,c,f} \P_p(0 \leftrightarrow b)\, C_3^{\sss (0)} (b,c,f,b)\, 
	\P_p(c \leftrightarrow f) \sum_{v: \{a,v\} \in \Ecal} p\, \P_p(f \leftrightarrow v)\nn\\
	&\qquad= \sup_{a'} \sum_{b} C_2^{\sss (1)}(0,b,a') \sum_{c,f} C_3^{\sss (0)} (b,c,f,b)\, 
	\P_p(c \leftrightarrow f),\nn
	}
where for the second equality we use translation invariance of the two-point function to shift the vertex 
$a$ to $a'= a-f+b$ and $f$ to $b$ in the term $\sum_{v: \{a,v\} \in \Ecal} p\, \P_p(f \leftrightarrow v)$. 
Using that $\P_p(c \leftrightarrow f) \le 1$ and applying Lemma~\ref{lem:poly} twice, we find that the 
remaining terms are bounded by $O(m^{-1} + V^{-1/3})$. Combined with the bounds on the other two 
terms, we thus conclude that the case where $h=k$ and $f=g$ contributes $O(m^{-2} + V^{-1/3})$.

It remains to bound the case $h=k$ and $d(f,g)=1$. We apply Lemmas~\ref{lem:poly} and \ref{lem:distone} 
to bound the two sums in $D_3$ by $O(m^{-1} + V^{-1/3})$ and write the remaining $D_2$ terms as
\begin{multline}
	\sup_a \sum_{b,c,f} \sum_{g: \{f,g\} \in \Ecal} C_5^{\sss (1)}(0,b,g,c,f,a)\,
	 \P_p(b \leftrightarrow c) \\
	\le \sup_{a'} \sum_{b,f'} C_3^{\sss (1)}(0,b,f',a') \sup_{f}  \sum_{c,g : \{f,g\} \in \Ecal} 
	C_3^{\sss (0)}(b, c,g,b),
\end{multline}
where for the bound we use the translation invariance of the two-point function again to shift $a$ to 
$a' = a-c + b$, $f$ to $f' = f -c + b$ and $c$ to $b$, and we take the supremum over $f$ in the second 
sum for an upper bound. Apply Lemma~\ref{lem:distone} to the second sum and Lemma~\ref{lem:poly} 
to the first, to bound this factor by $O(m^{-1} + V^{-1/3})$ also. 

We have thus bounded all contributions to $F_2$, and conclude that
\begin{equation}\label{e:F2bd}
	F_2  = O(m^{-2} + V^{-1/3}).
\end{equation}

The remaining terms in \eqref{e:PiNdiag} are $O(\theta^{N-2})$. This is proved by Borgs et al.\ 
\cite{BorChaHofSlaSpe05b}. More precisely, combining \cite[(4.42)]{BorChaHofSlaSpe05b} with \eqref{e:triangledef}, \eqref{e:pcdef}, \eqref{e:opentri}, \eqref{e:F1bd} and \eqref{e:F2bd}, we may 
conclude that for all $N \ge 2$,
\begin{equation}
	\hat \Pi_p^{\sss (N)} = (1 + 10 \theta^3) \big(3 \theta^3 (1+ 10 \, \theta^3)\big)^{N-2} 
	O(m^{-2} + V^{-1/3}).
\end{equation}
Hence, the sum $\sum_{N=2}^\infty \hat \Pi_p^{\sss (N)} $ is of order $O(m^{-2} + V^{-1/3})$ when 
$3 \theta^3 (1+ 10 \, \theta^3) < 1$, which is the constraint on $\theta$ that we assumed. \qed

\section*{Acknowledgments}

The work in this paper is supported by the Netherlands Organisation for Scientific Research 
(NWO) through Gravitation-grant NETWORKS-024.002.003. RvdH is also supported by 
NWO through VICI-grant 639.033.806, and FdH by the European Research Council (ERC) 
through Advanced Grant VARIS-267356.


\begin{small}
\bibliographystyle{abbrv}
\bibliography{..//LorenzosBib}
\end{small}
\end{document}